\documentclass{amsproc}
\usepackage{color}
\newtheorem{Theorem}{Theorem}
\newtheorem{theorem}{Theorem}[section]
\newtheorem{lemma}[theorem]{Lemma}
\newtheorem{proposition}[theorem]{Proposition}
\newtheorem{corollary}[theorem]{Corollary}

\theoremstyle{definition}
\newtheorem{definition}[theorem]{Definition}
\newtheorem{example}[theorem]{Example}

\newtheorem{assumption}[theorem]{Assumption}

\newtheorem{question}[theorem]{Question}

\theoremstyle{remark}
\newtheorem{remark}[theorem]{Remark}

\numberwithin{equation}{section}

\newcommand{\abs}[1]{\lvert#1\rvert}

\newcommand{\bd}{ {\sc Proof}.\ \ }

\begin{document}

\title[]{Dynamics of isolated left orders}


\author{Shigenori Matsumoto}
\address{Department of Mathematics, College of
Science and Technology, Nihon University, 1-8-14 Kanda-Surugadai,
Chiyoda-ku, Tokyo, 101-8308 Japan}
\email{matsumo@math.cst.nihon-u.ac.jp}
\thanks{2010 {\em Mathematics Subject Classification}. Primary 20F65.
secondary 20F05.}
\thanks{{\em Key words and phrases.}  left orderable groups, isolated
left orders}

\thanks{The author is partially supported by {{Grant-in-Aid for
Scientific Research}} (C) No.\ 25400096.}
%    General info

\date{\today}

\newcommand{\AAA}{{\mathcal A}}
\newcommand{\BBB}{{\mathbb B}}
\newcommand{\LL}{{\mathcal L}}
\newcommand{\MCG}{{\rm MCG}}
\newcommand{\PSL}{{\rm PSL}}
\newcommand{\R}{{\mathbb R}}
\newcommand{\Z}{{\mathbb Z}}
\newcommand{\XX}{{\mathcal X}}
\newcommand{\per}{{\rm per}}
\newcommand{\N}{{\mathbb N}}

\newcommand{\PP}{{\mathcal P}}
\newcommand{\RR}{{\mathcal R}}
\newcommand{\GG}{{\mathbb G}}
\newcommand{\FF}{{\mathcal F}}
\newcommand{\EE}{{\mathbb E}}
\newcommand{\BB}{{\mathbb B}}
\newcommand{\CC}{{\mathcal C}}
\newcommand{\HH}{{\mathcal H}}
\newcommand{\UU}{{\mathcal U}}
\newcommand{\oboundary}{{\mathbb S}^1_\infty}
\newcommand{\Q}{{\mathbb Q}}
\newcommand{\DD}{{\mathcal D}}
\newcommand{\rot}{{\rm rot}}
\newcommand{\Cl}{{\rm Cl}}
\newcommand{\Index}{{\rm Index}}
\newcommand{\Int}{{\rm Int}}\newcommand{\Fix}{{\rm Fix}}
\newcommand{\Fr}{{\rm Fr}}
\newcommand{\ZZ}{\Z[2^{-1}]}
\newcommand{\II}{{\mathcal I}}
\newcommand{\JJ}{{\mathcal J}}
\newcommand{\g}{{\rm gen}}
\newcommand{\KK}{{\mathcal K}}
\newcommand{\OO}{{\mathcal O}}
\newcommand{\MM}{{\mathcal M}}
\newcommand{\NN}{{\mathcal N}}
\newcommand{\oM}{{\rho_M}}
\newcommand{\oG}{{\overline G}}
\newcommand{\oH}{{\overline {\mathcal H}}}
\newcommand{\og}{{\overline g}}

\date{\today }

\maketitle

\begin{abstract}
 A left order of a countable group $G$ is called isolated if it is an
 isolated point in the compact space $LO(G)$ of all the left orders of
 $G$. We study properties of a dynamical realization of an isolated left
 order. Especially we show that it acts on $\R$ cocompactly.
As an application, we give a dynamical proof of the Tararin theorem
 which characterizes those countable groups which admit only finitely many
left orders. We also show that the braid group $B_3$ admits countably
 many isolated left orders which are not the automorphic images of the 
others.
\end{abstract}

\section{Introduction}
Throughout this paper, all the groups considered are countable.
Given a group $G$, a total order $<_\lambda$ on $G$ is called a {\em left order}
if for any $f,g,h\in G$, $f<_\lambda g$ implies $hf<_\lambda hg$. An
element $g\in G$ is called $\lambda$-positive if $g>_\lambda e$.
The set of all the $\lambda$-positive elements is called the {\em positive cone of
$\lambda$} and is denoted by $P_\lambda$. It is a subsemigroup and 
$P_\lambda\sqcup P_\lambda^{-1}=G\setminus\{e\}$.

Given a left order $<_\lambda$, we define $\lambda:G\setminus\{e\}\to\{\pm1\}$ by
$\lambda(g)=1$ if and only if $g\in P_\lambda$.
 Then we  have
\begin{equation}\label{e1}
\lambda(f)=1, \lambda(g)=1 \Rightarrow \lambda(fg)=1, \mbox{ and }
\lambda(f^{-1})=-\lambda(f).
\end{equation}
Conversely given a map $\lambda:G\setminus\{e\}\to\{\pm1\}$ which
satisfies (\ref{e1}), we get a left order $<_\lambda$ by setting
$f<_\lambda g$ if $\lambda(f^{-1}g)=1$.
The map $\lambda$ is also referred to as a left order.
Thus the set $LO(G)$ of the left orders on $G$ is viewed as a closed subset of
the space $\{\pm1\}^{G\setminus\{e\}}$ with the pointwise convergence topology.
This yields a totally disconnected compact metrizable topology on
$LO(G)$
(metrizable since $G$ is countable).
It is either finite or uncountably many \cite{L}. We call $\lambda\in
LO(G)$ {\em isolated} if it is an isolated point in the space $LO(G)$.

Given $\lambda\in LO(G)$, there is defined a dynamical realization 
$$\rho_\lambda:G\to {\rm Homeo}_+(\R)$$ 
 based at $x_0\in\R$ such that $f<_\lambda g$ if and only if
 $fx_0<gx_0$.
We discuss its fundamental properties in Section 2.
Especially we show that the dynamical realization is tight at the base
point. See Definition \ref{def-tight}.

In this paper, we are mainly interested in isolated orders, since in
this case, the dynamical realizations display a certain kind of rigidity,
and vice versa. {In \cite{MR} Theorems 1.2 and 3.11, the relation between
the isolation of left orders and the rigidity of the dynamical realization is
described, as well as for circular orders.}

An action $\rho:G\to {\rm Homeo}_+(\R)$ is said to be {\em cocompact} if
there is a compact interval $I$ such that any orbit $\rho(G)x$
intersects $I$. Our first result, proved in Section 3, is the following.

\begin{Theorem}\label{t1}
 If $\lambda\in LO(G)$ is isolated, then its dynamical realization
 $\rho_\lambda$ is cocompact.
\end{Theorem}

{ In fact if the group $G$ is finitely generated,
the dynamical realization is cocompact for any left order, isolated or
not.
(See Lemma 3.1 below.)
 Therefore Theorem \ref{t1} is mainly concerned with non finitely
 generated groups.}
By Theorem \ref{t1}, the dynamical realization of an isolated order
admits a minimal set $\MM$, which is shown to be unique unless $G\cong\Z$.
In Section 4, we show that if $\MM=\R$, then the group is rational
(Theorem \ref{t31}).

Given $\lambda\in LO(G)$, a subgroup $H$ of $G$ is called $\lambda$-{\em convex},
if whenever $h_1,h_2\in H$, $g\in G$ and $h_1<_\lambda g <_\lambda h_2$, we have
$g\in H$. The set of convex subgroups is totally ordered by the
inclusion. The following theorem is shown in Section 5.

\begin{Theorem}\label{t2}
 If $\lambda\in LO(G)$ is isolated, then there are only finitely many 
$\lambda$-convex subgroups.
\end{Theorem}

{This is known to specialists (see for example
\cite{DNR} Exercise 3.3.15). However, our strategy of the proof is
different from that mentioned in \cite{DNR}.}

Theorem \ref{t2} enables us to define the maximal sequence of convex subgroups of an
isolated left order. 
As an application of our method, we
give a dynamical proof of the Tararin theorem which characterizes
the groups with finitely many left orders in Section 6.
In Section 7, the maximal Tararin subgroup of an isolated left order is
defined, and is shown to be equal to the Conradian soul \cite{N}.

Last sections 8 and 9 are more or less independent of the previous sections.
{ 
Dubrovina-Dubrovin \cite{DD} constructed an isolated order
 $\lambda_n$ on the braid group $B_n$, $n\geq3$.
In section 9, we show:

\begin{Theorem}\label{t3}
 There are countably many
 isolated orders in $LO(B_3)$ which
are not the automorphic images
\footnote{Given $\phi\in{\rm Aut}(G)$ and $\lambda\in LO(G)$, the left
 order $\phi^*\lambda\in LO(G)$ defined by 
$g<_{\phi^*\lambda} g'$
 if and only if $\phi(g)<_\lambda \phi(g')$
is called an automorphic image of
 $\lambda$.
 For example, the reciprocal of the natural order $\lambda$ of $\Z$ is
 an automorphic image of $\lambda$.}
 of the others.
\end{Theorem}

The method is a modification of the proof of \cite{MR} Theorem 1.4.
The following theorem is the starting point of the proof of
Theorem \ref{t3}.
Let
$$
G=\langle a,b\mid a^2=b^3\rangle,\ \ \
\overline G=\langle\alpha,\beta\mid\alpha^2=\beta^3=e\rangle,$$
and $q:G\to\overline G$ the surjective homomorphism defined by
$q(a)=\alpha$ and $q(b)=\beta$. Notice that $G$ is isomorphic to $B_3$, and $\oG$ to $PSL(2,\Z)$.
We denote by $CO(\oG)$ the space of the left invariant circular orders of
$\oG$ (see Section 8). In Section 9, we show:

\begin{Theorem}\label{t4}
 The homomorphism $q$ induces a homeomorphism $q_*:LO(G)\to CO(\oG)$.
\end{Theorem}}

Isolated left orders are often induced from isolated circular orders of
the group quotiented by the center. 
See \cite{MR} Section 5, for example. The above theorem is also a typical example.
However there is an example of a group
with isolated left orders which admits no center,  
constructed in \cite{I2} 3.2.

\medskip
{ 
{\sc Acknowledgement.} The author is grateful to Y. Matsuda for stimulating
conversations. Hearty thanks are due to the referee for many helpful
suggestions.
}

\section{Dynamical Realization}

In this section, we define a dynamical realization of a left order
$\lambda\in LO(G)$ and study its fundamental properties.
Fix an enumeration of $G$: $G=\{g_i\mid i\in\N\}$ such that $g_1=e$.
We define an order preserving embedding $\iota:G\to\R$ inductively as
follows. Define $\iota(g_1)=x_0$, where $x_0$ is some point in $\R$. Assume we have defined $\iota$ on the subset
$\{g_1,\ldots,g_n\}$, $n\geq1$, and let us define $\iota(g_{n+1})$.
Order the subset $\{g_1,\ldots,g_n\}$ as
$$ 
g_{i_1}<_\lambda g_{i_2}<_\lambda\cdots<_\lambda g_{i_n}.$$

If $g_{n+1}<_\lambda g_{i_1}$, define $\iota(g_{n+1})=\iota(g_{i_1})-1$, 

if $g_{i_n}<_\lambda g_{n+1}$,  $\iota(g_{n+1})=\iota(g_{i_n})+1$,

and if $g_{i_k}<_\lambda g_{n+1}<_\lambda g_{i_{k+1}}$, 
$\iota(g_{n+1})=(1/2)(\iota(g_{i_k})+\iota(g_{i_{k+1}}))$.

Then we have $\inf\iota(G)=-\infty$ and $\sup\iota(G)=\infty$.
The left translation of $G$ yields an order preserving action of $G$ on
$\iota(G)$, which extends to a continuous action on the closure 
$\Cl(\iota(G))$. (See Proposition 2.2 and the proof of Corollary 2.3.)
Extend it further to a continuous action on $\R$ by setting
that the action on gaps of $\Cl(\iota(G))$ be linear. 

This action is called the {\em dynamical realization of} $\lambda$ {\em
based at $x_0$}, and is denoted by $\rho_\lambda$.
The dynamical realization depends on the choice of the enumeration of $G$.
 Soon later, we
shall show that any two dynamical realizations are mutually topologically conjugate. 

\begin{definition}
 \label{def-tight}
An action $\rho:G\to {\rm Homeo}_+(\R)$ is called {\em tight at}
 $x_0\in\R$ if

(1) $\rho$ is free at $x_0$ i.e, the stabilizer at $x_0$ is trivial,

(2) $\inf\rho(G)x_0=-\infty$, $\sup\rho(G)x_0=\infty$, and

(3) whenever $\Cl(\rho(G)x_0)\cap[a,b]=\{a,b\}$ for any $a<b$, we have
$\{a,b\}\subset\rho(G)x_0$.
\end{definition}

\begin{proposition}\label{prop-tight}
The dynamical realization $\rho_\lambda$ based at $x_0$ is tight at $x_0$.
 \end{proposition}

\bd All that needs proof is (3). Let $a<b$ be as in (3). The proof is by
contradiction. Assume, to fix the idea, that
$a\not\in\rho(G_0)x_0=\iota(G)$.
(Notice that $\iota(g)=\rho_\lambda(g)x_0$.)
Choose $\epsilon$ small enough compared with $b-a$, and choose
$\iota(g_1)\in (a-\epsilon,a)$ and $\iota(g_2)\in[b,b+\epsilon)$.
Recall that the dynamical realization is defined via an enumeration of
$G$. One may assume that there is no point in
$(\iota(g_1),\iota(g_2))\cap\iota(G)$ which is enumerated before $g_1$
or $g_2$, since otherwise one may pass to that point. Since
$a\in\Cl(\iota(G))\setminus\iota(G)$, there is a point $\iota(g_3)$ in 
$(\iota(g_1),\iota(g_2))\cap\iota(G)$ which is enumerated for the first
time after $g_1$ and $g_2$. Then $\iota(g_3)$ is the midpoint of
$\iota(g_1)$ and $\iota(g_2)$ and must be fallen in $(a,b)$ since
$\epsilon$ is small. A contradiction. \qed

\begin{corollary}\label{c25}
The dynamical realizations defined via two different enumerations of $G$
are mutually conjugate by an orientation and base point preserving
 homeomorphism of $\R$ 
\end{corollary}

\bd Let $\iota$ and $\iota'$ be two embeddings of $G$ obtained by
different enumerations of $G$. There is an orientation preserving
bijection $h:\iota(G)\to\iota'(G)$ defined by $h(\iota(g))=\iota'(g)$
($g\in G$). By the tightness, $h$ extends, first of all. to a homeomorphism
$h:\Cl(\iota(G))\to\Cl(\iota'(G))$, and then to a homeomorphism of $\R$ linearly
on gaps. The extended $h$ yields the required conjugacy. \qed

\medskip
 The proof of the previous corollary also yields the
following result, which will be used in Section 9.

\begin{corollary}
 \label{bijection}
Let $\HH$ be the set of the orientation and base point $x_0$ preserving
topological conjugacy classes of the homomorphisms 
$G\to{\rm  Homeo}_+(\R)$
which are tight at $x_0$. Then the dynamical realization at $x_0$ induces a
bijection of $LO(G)$ onto $\HH$. \qed
\end{corollary}

\bigskip
A left order $<_\lambda$ is called {\em discrete} if there is a minimal 
$\lambda$-positive element, and {\em indiscrete} otherwise. 

\begin{corollary}
 \label{c21}
If $\lambda\in LO(G)$ is indiscrete, then the orbit $\rho_\lambda(G)x_0$
 of the base point $x_0$ is dense in $\R$.
\end{corollary}

\bd Assume $\Cl(\rho_\lambda(G)x_0)\neq\R$ and let $(a,b)$ be a gap of
$\Cl(\rho_\lambda(G)x_0)$. Then by the previous lemma, we have 
$a,b\in\rho_\lambda(G)x_0$. That is, $a=\iota(g_1)$  and $b=\iota(g_2)$
for some $g_1,g_2\in G$.
Then $g_1^{-1}g_2$ is the minimal positive element, and $\lambda$ is
discrete.  \qed

\section{Proof of Theorem \ref{t1}}

We begin with two lemmas. { The first one can be found in \cite{N2}
(Proposition 2.1.12).} 

\begin{lemma}\label{l21}
 Let $G$ be a finitely generated group which acts on $\R$ without global
 fixed points. Then the action is cocompact.
\end{lemma}

\bd 
We identify $\R\approx (0,1)$. Let $G_0$ be a finite generating set of
$G$. Define 
$$ a=\sup_{s\in G_0}\sup_{x\in(0,1)}\vert sx-x\vert.
$$
Choose a compact interval $J\subset(0,1)$ such that $\vert J\vert>a$.
Given any point $x\in(0,1)$, we have $\inf Gx=0$ and $\sup Gx=1$ since
there is no global fixed point. Considering the Schreier graph of $Gx$,
one can show that $Gx\cap J\neq\emptyset$. \qed

\begin{lemma}\label{l22}
Let $G$ be a group acting on $\R$ and let $y_0\in \R$. Denote by $G_{y_0}$
  the stabilizer of $G$ at $y_0$. { Assume $G_{y_0}\neq G$.}
Given $\lambda_0\in LO(G_{y_0})$,
 there are at least two orders in $LO(G)$ which restrict to
 $\lambda_0$ on $G_{y_0}$.
\end{lemma}

\bd Let $\mu$ be the $G$-invariant order on $G/G_{y_0}$ given by the natural
order of the orbit $Gy_0\approx G/G_{y_0}$ in $\R$. Then $\lambda_0$ and $\mu$
determines a left order on $G$ lexicographically (Lemma \ref{l41}). If we consider the reciprocal order $-\mu$,
we get another one. \qed

\bigskip

Assume $\lambda\in LO(G)$ is an isolated left order on $G$. Since we are considering the
pointwise convergence topology, this is equivalent to the following
condition ($\star$)

\smallskip
($\star$)\em There is a finite subset $S\subset P_\lambda$ such that
$\lambda$ is the only element in $LO(G)$ which contains $S$ in its
positive cone. 

\rm\smallskip
Such a subset $S$ is called a {\em characteristic positive set} of $\lambda$.

\bigskip
\sc Proof of Theorem \ref{t1}.\ \rm
By the dynamical realization of the isolated left order $\lambda$, the group $G$ acts on $\R$.
Let $H$ be the subgroup of $G$ generated by a characteristic
positive set $S$ of $\lambda$. If there is no global
fixed point by the action of $H$, then $H$ acts on $\R$ cocompactly by
Lemma \ref{l21}, and
hence also $G$, finishing the proof. In the remaining case, choose a
global fixed point $y_0$ of $H$ and consider $G_{y_0}$.
{  We have $G_{y_0}\neq G$ since the dynamical realization has
no global fixed point, by its tightness.} By the previous lemma,
the restriction of $\lambda$ to $G_{y_0}$ extends
to two left orders of $G$. 
But we have $S\subset H\subset G_{y_0}$ and hence $S$ is contained in
the positive cone of both orders.
A contradiction. \qed

\begin{remark}
The condition that $\lambda$ be isolated is actually necessary for
 Theorem \ref{t1}. To show this, let $G$ be the infinite direct sum of
 $\Z$, i.e,
$$
G=\{(a_n)_{n\in\N}\mid a_n\in\Z,\ \ a_n=0\ \mbox{ but for finitely many }n\}.$$
Define a left order on $G$ by setting $0<(a_n)$ if $0<a_N$, where
$N$ is the largest number such that $a_N\neq 0$.
Then its dynamical realization is not cocompact. To show this, define for
 $m\in\N$, 
$$
G_m=\{(a_n)\mid a_n=0,\ \forall n>m\}.$$
Then $G_m$'s form an exhausting increasing sequence of convex subgroups. 
Consider the dynamical realization $\rho_\lambda$ based at $x_0$.
The points 
$$
\xi_n=\inf\rho_\lambda(G_n)x_0\ \mbox{ and }\ 
\eta_n=\sup\rho_\lambda(G_n)x_0
$$
are fixed points of $\rho_\lambda(G_n)$. They satisfy
 $\xi_n\searrow-\infty$ and $\eta_n\nearrow\infty$ by condition (2) of
Definition \ref{def-tight}, since $G_n$ is exhausting.
This implies that $\rho_\lambda$ is not cocompact.
\end{remark}

\bigskip

Theorem \ref{t1} implies that there is a minimal set $\MM$ for the dynamical
realization of an isolated left order.\footnote{ 
{  The results of the
 remaining part of this section (trichotomy, Lemma 3.4 and
Corollary 3.5) and Lemma 5.4 (1) (2) hold true whenever the left order
$\lambda$ admits a cocompact dynamical realization, especially when
$G$ is finitely generated. But we shall state it only for an isolated
left order $\lambda$.}}
There is a trichotomy for $\MM$ (\cite{DKNP} Proposition 6.1).

\smallskip
(I) \em  $\MM=\R$. \rm

(II) \em $\MM$ is infinite and discrete in $\R$.\rm

(III) \em $\MM$ is locally Cantor. In this case, if $X$ is a nonempty closed
subset of $\R$ invariant by the dynamical realization of $G$, then
$\MM\subset X$. Especially, $\MM$ is the unique minimal set.\rm

\begin{lemma}
 \label{l35}
Let $\lambda\in LO(G)$ be isolated, with $\MM$ an associated minimal set.
Assume (III) above, or (II) and $G\not\cong\Z$. Then the base point
 $x_0$ is contained in a gap $I_1$ of $\MM$, the stabilizer $G_{I_1}$ is
nontrivial, and there is no gap of $\MM$ other than the orbit of $I_1$.
\end{lemma}

\bd We give a proof only for case (III). { Case (II) can be
treated much easier. Notice that if $\rho(G)x_0$ is discrete,
then $G\cong\Z$.}
Assume that the base point $x_0$ is contained in $\MM$ and let $(a,b)$
be a gap of $\MM$. Since the dynamical realization $\rho_\lambda$ is
tight, we have $a,b\in\rho_\lambda(G)x_0$.
But there is no orientation preserving homeomorphism leaving $\MM$
invariant and mapping $a$ to $b$. The contradiction shows that $x_0$
is contained in a gap $I_1$ of $\MM$.

If $G_{I_1}$ is trivial, then $\rho_\lambda(G)x_0\cap I_1=\{x_0\}$.
Again by the tightness, the boundary points of $I_1$ must belong to
$\rho_\lambda(G)x_0$. A contradiction. The last statement follows
similarly from the tightness. \qed

\begin{corollary}\label{minimal-unique}
If $G\not\cong\Z$, then the minimal set $\MM$ of the dynamical
 realization $\rho_\lambda$ of an isolated left order $\lambda$ is unique.
\end{corollary}

\bd
 { All that needs proof is the case where $\MM$ is discrete,
 since a locally Cantor minimal set is always unique.
We still use the notation of the previous lemma. 
If there is another minimal set $\MM'$,  then $\MM'\cap I_1$ must be one
point, say $y_0$, which is fixed by $\rho_\lambda(G_{I_1})$.
But then $\rho_\lambda(G)x_0\cap I_1=\rho_\lambda(G_{I_1})x_0$ must be contained in an open subinterval of 
$I_1$ delimited by $y_0$, contrary to the tightness.
}
 \qed

\section{The case $\MM=\R$}

This section is devoted to the proof of the following theorem.

\begin{theorem}\label{t31}
Let $\lambda\in LO(G)$ be isolated and assume that
the dynamical realization $\rho_\lambda$ is minimal. Then
the group $G$ is isomorphic to an additive subgroup $A$ of $\Q$ such
 that
$A\not\cong\Z$,
{ and  $\lambda$ is either the natural left order given by 
$A\subset\Q\subset\R$
or its reciprocal.
}
\end{theorem}

%First of all, we shall show the following. Let $A$ be a dense subgroup of $\Q$.
%\begin{lemma}\label{l31}
% A minimal action of $A$ on $\R$ is conjugate either to the translation
% by $A$ or by $-A$.
%\end{lemma}

%\bd Choose any $a_0\in A\setminus\{0\}$. The action of the group 
%$\langle a_0\rangle$ is free. For, if not, the fixed point set
%$\Fix( a_0)$ would be invariant by the whole group $A$, 
%since $A$ is abelian. Consider the action of the group 
%$A/\langle a_0\rangle$ on $\R/\langle a_0\rangle\approx S^1$. Since
%$A/\langle a_0\rangle$ is amenable, there is an 
%$A/\langle a_0\rangle$-invariant probability measure. Its lift $\mu$
%is a locally finite $A$-invariant measure. Since the action of $A$ is
%minimal, $\mu$ is atomless and fully supported. Thus there is a
%homeomorphism $h$ such that $h_*\mu$ is the Lebesgue. 
%Conjugating the $A$-action by $h$, we get an action by translations,
%Now any nontrivial  homomorphism from $A$ to $\R$ is the
%multiplication by a nonzero constant. Conjugating by an appropriate linear map,
%we get the translation by $A$ or $-A$. \qed

{  This theorem might be known among specialists, but the
author cannot locate it in the literature.}

Let $\lambda$ be an element of $LO(G)$ which satisfies the hypothesis of Theorem \ref{t31}. 
We shall abbreviate the notations $\rho_\lambda(g)x$ by $gx$,
and $\rho_\lambda(G)\subset{\rm Homeo}_+(\R)$ by $G$.
Let $Z$ be the centralizer of $G$ in ${\rm Homeo}_+(\R)$.

\begin{lemma}\label{l32}
The centralizer $Z$ is an abelian group which acts freely on $\R$.
\end{lemma}

\bd For $\zeta\in Z\setminus\{\rm id\}$, $\Fix(\zeta)$ is a closed set
which is invariant by $G$. Since the $G$-action is minimal, we have
$\Fix(\zeta)=\emptyset$. By H\"older's theorem (e.g, \cite{N}), any group acting freely
on $\R$ is abelian. \qed

\bigskip
 Let $x_0$ be the base point of the dynamical realization.
Choose $x_n\in Gx_0$, $n\in\N$, such that $x_n\to x_0$, $x_n\neq x_0$.
Notice that $G$ acts freely at $x_n$.
Let $\lambda_n\in LO(G)$ be the order determined by $x_n$:
$g>_{\lambda_n}e$ if and only if $gx_n>x_n$.
Then $\lambda_n\to\lambda$ in $LO(G)$.
Since $\lambda$ is isolated, $\lambda_n=\lambda$ for any large $n$. We
assume $\lambda_n=\lambda$ for all $n$. 
Define an order preserving bijection $\zeta_n: Gx_0\to Gx_n$ by
$\zeta_n(gx_0)=gx_n$. Since $Gx_0=Gx_n$ is dense in $\R$, the map
$\zeta_n$ extends to an orientation preserving homeomorphism of $\R$, denoted by
the same letter $\zeta_n$. Clearly $\zeta_n\neq id$.

\begin{lemma}\label{l33}
 We  have $\zeta_n\in Z$.
\end{lemma}

\bd Given any $g\in G$, it suffices to show that $\zeta_ng=g\zeta_n$ on the dense
subset $Gx_0$. For any $hx_0\in Gx_0$, we have
$$
\zeta_ng(hx_0)=\zeta_n((gh)x_0)=ghx_n=g(hx_n)=g\zeta_n(hx_0),
$$
as is required. \qed

\begin{lemma}\label{l34}
The action of $Z$ is minimal, and is conjugate to translations. 
\end{lemma}

\bd By Lemmas \ref{l32} and \ref{l33}, there is an element in $Z$ which
acts freely on $\R$. This implies that the action of $Z$ is cocompact.
Let $\NN$ be a minimal set of $Z$. If it is locally Cantor, then $\NN$
is the unique minimal set, and must be invariant by $G$. But $G$-action is
minimal
by the assumption. A contradiction. Next assume $\NN$ is discrete.
Then since the $Z$-action is free, we must have $Z\cong\Z$,
contradicting Lemma \ref{l33}.
Therefore $Z$ must act minimally on $\R$. 

Choose any $\zeta_0\in Z\setminus\{\rm id\}$. Since the action of the group 
$\langle \zeta_0\rangle$ is free and
 $Z$ is abelian,  the group 
$Z/\langle \zeta_0\rangle$ acts on $\R/\langle \zeta_0\rangle\approx
S^1$.
Since $Z/\langle\zeta_0\rangle$ is amenable, there is an 
$Z/\langle \zeta_0\rangle$-invariant probability measure. It lifts
to a locally finite $Z$-invariant measure $\mu$ on $\R$. Since the action of $Z$ is
minimal, $\mu$ is atomless and fully supported. Thus there is a
homeomorphism $h$ such that $h_*\mu$ is the Lebesgue. 
Conjugating the $Z$-action by $h$, we obtain an action by translations.
\qed

\bigskip
\sc Proof of Theorem \ref{t31}.
\rm  By changing the coordinate, we assume that the action of $Z$ is by
translations. Since the $Z$-action is minimal, any element of
$G$, commuting with $Z$, acts also by
translations. Then we have an injective homomorphism $\phi:G\to\R$ defined by
the translation length. We shall show that $\phi$ embeds $G$ into $\Q$.
Assume not. Then $G$ is a nontrivial direct sum: $G=G_1\oplus G_2$.
Given any $a\in\R$, we obtain a homomorphism $\phi_a:G\to\R$ by
setting $\phi_a=\phi$ on $G_1$ and $\phi_a=a\phi$ on $G_2$. There is $a$
arbitrarily near 1 such that $\phi_a$ is injective. But $\phi_a$ yields
a left order different from $\lambda$ and arbitrarily near $\lambda$. This contradicts the
assumption that $\lambda$ is isolated,
finishing the proof that $G$ is isomorphic to an additive subgroup $A$ of
$\Q$. { The last statement of the theorem follows at once. 
}\qed

\section{Convex subgroups}

We shall prove Theorem \ref{t2} in this section. First we begin
with fundamental properties of convex subgroups. For the definition
of convex subgroups, see Introduction.
We begin with a well known easy fact.
\begin{lemma}\label{l41}
 Let $H$ be a subgroup of $G$. For any $\lambda_0\in LO(H)$
and any $G$-invariant total order $\lambda_1$ on $G/H$, there is
a unique order $\lambda\in LO(G)$ such that $H$ is $\lambda$-convex,
that $\lambda\vert_H=\lambda_0$,
and that for $g\not\in H$, $g>_{\lambda}e$ if and only if $gH>_{\lambda_1}H$.
\qed
\end{lemma}

Such an order $\lambda$ is said to be {\em determined lexicographically
by $\lambda_0$ and
$\lambda_1$.}

\begin{lemma}\label{l42}
Let $\lambda\in LO(G)$ and $H$ a $\lambda$-convex subgroup of $G$. 
Then there is a $G$-invariant total order $\lambda_1$ on $G/H$ such that
$\lambda$ is determined lexicographically by
$\lambda\vert_H$ and $\lambda_1$.
 \end{lemma}

\bd 
Define  a total order $\lambda_1$ on $G/H$ by setting $g_1H<_{\lambda_1}g_2H$ if
$e<_\lambda g_1^{-1} g_2$ and $g_1^{-1}g_2\not\in H$. The convexity of $H$ shows that this is a well defined
$G$-invariant order. \qed

\bigskip
{ If $G$ is isomorphic to $\Z$ or if the minimal set of the
dynamical realization of $\lambda$ is $\R$, then there is no proper
$\lambda$-convex subgroups, and Theorem \ref{t2} holds true.}
Henceforth in this section we work under the following assumption.

\begin{assumption}
(1) $\lambda\in LO(G)$ is isolated with
a characteristic positive set $S$.  

(2) $G$ is not isomorphic to $\Z$.

(3) The minimal set $\MM$ of the dynamical realization is not $\R$.
\end{assumption}

 Denote by $I_1=(y_0,z_0)$ the gap of $\MM$ which contains the base
 point $x_0$ (Lemma  \ref{l35}), and by $G_1$ the stabilizer of $I_1$.

\begin{lemma}\label{l43}
(1) $G_1$ is proper and nontrivial.

 (2) $G_1$ is the maximal proper $\lambda$-convex subgroup of $G$.

(3) The restricted order $\lambda\vert_{G_1}$ is isolated with
 characteristic positive set $S\cap G_1$.

(4)  $S\cap(G\setminus G_1)\ne\emptyset$.
\end{lemma}

\bd The subgroup $G_1$ is clearly proper. It is nontrivial by Lemma \ref{l35}.
Also $G_1$ is convex. Let $H$ be an arbitrary
proper $\lambda$-convex subgroup of $G$. We shall show that 
$H\subset G_1$. Consider first the case where $\MM$ is discrete.
By looking at the action of $G$ on $\MM$, one can define a surjective
homomorphism $\phi:G\to\Z$ such that ${\rm Ker}(\phi)=G_1$. 
If $\phi(H)$ is nontrivial, then clearly we have $H=G$ since $H$ is
convex. If $\phi(H)$ is trivial, then $H\subset G_1$, as is required.

So in the rest, we assume that $\MM$ is locally Cantor.
Let $\HH$ be the convex hull of $Hx_0$ in $\R$.
{  Then $\HH$ is a bounded open interval of $\R$.
The boundedness follows from the convexity and the properness of $H$}. The convexity of $H$ implies that
for any $g\in G$, we have either
$g\HH=\HH$ or $g\HH\cap\HH=\emptyset$. Thus the closed set
$$X=\R\setminus \bigcup_{g\in G}g\HH$$ is $G$-invariant and nonempty.
Therefore we have $\MM\subset X$, which implies $\HH \subset I_1$,
showing that $H\subset G_1$.

Let us show that  $S\cap G_1$ is a characteristic positive set of
$\lambda\vert_{G_1}$. If not, there is a left order $\lambda'_0$
($\lambda'_0\neq\lambda\vert_{G_1}$) of $G_1$
such that $S\cap G_1$ is contained in the positive cone of $\lambda'_0$.
Let $\lambda_1$ be the $G$-invariant total order on $G/G_1$ obtained by
Lemma \ref{l42}. Let 
$\lambda'\in LO(G)$ be the order determined lexicographically by
$\lambda'_0$ and $\lambda_1$. Then $\lambda'$ contains $S$ in its
positive cone and $\lambda'\neq\lambda$, contradicting that $S$ is a
characteristic positive set of $\lambda$.

Finally let us show that $S\cap(G\setminus G_1)$ is nonempty.
If it is empty, then $\lambda\vert_{G_1}$ and $-\lambda_1$ lexicographically
determines $\lambda'\in LO(G)$, where $-\lambda_1$ is the reciprocal of
the order $\lambda_1$ constructed in Lemma \ref{l42}. 
But $S$ is contained in the positive cone of $\lambda'$. A contradiction.
\qed

\bigskip

\sc Proof of Theorem 2. \rm
By Lemma \ref{l43}, we obtain the maximal proper convex subgroup $G_1$.
If $G_1$ is not isomorphic to $\Z$ and the minimal set of
the dynamical realization of 
$\lambda\vert_{G_1}$ is not the whole $\R$, then we can repeat the
process and obtain the second maximal proper convex subgroup $G_2$.
This process ends at finite steps since each time the number of elements
of positive characteristic set decreases. \qed

\begin{definition}
The sequence 
$$
G=G_0>G_1>\cdots >G_n>\{e\}$$
 of all the $\lambda$-convex subgroups is called the {\em maximal convex
 sequence} of the isolated order $\lambda$.
The number $n$ is called the {\em height} of $\lambda$.
\end{definition}

Thus an isolated left order with minimal dynamical realization has height 0.
Let $\MM_0$ be the minimal set of $G$ and $I_1$ the gap 
of $\MM_0$ containing the base point $x_0$. Then the maximal proper 
$\lambda$-convex
subgroup
$G_1$ is
the 
stabilizer of $I_1$. 
{ By Lemma 5.4 (3), $\lambda\vert_{G_1}$ is isolated,
and there is a minimal set $\MM_1$ of the $G_1$-action on
$I_1$.
}
Next consider the gap  $I_2$ of $\MM_1$ in $I_1$ containing $x_0$.
Continuing this way, we get 
a decreasing sequence of open intervals
$$
\R\supset I_1\supset\cdots\supset I_n.$$
Each subgroup $G_i$ is the stabilizer of $I_i$,
and each $\MM_i$ is a minimal set of $G_i$ in $I_i$.
The pair $(I_i,\MM_i)$ is called the {\em $i$-th internal pair
associated with the maximal convex sequence}.
There are only two
possibilities for the last group $G_n$:

(A)  $\MM_n=I_n$,

(B)  $G_n=\Z$.
\\
In (A), the order $\lambda$ is indiscrete and in (B), it is discrete.

{  As a corollary of Theorem \ref{t31}, we get the following
proposition, which will be used in the next section.

\begin{proposition}\label{new}
 If an isolated order $\lambda$ has height $0$, i.e, if there is no
 proper $\lambda$-convex subgroup, then the group $G$ is rational
and
the order $\lambda$ is the natural order of $G\subset\Q\subset\R$ or its reciprocal.
\end{proposition}}

\section{Tararin groups}

\begin{definition}
 A group $G$ is called a {\em Tararin group} if $\abs{LO(G)}<\infty$.
\end{definition}

Of course any left order of a Tararin group is isolated. In this section,
we shall give a dynamical proof of the following theorem by Tararin
\cite{T1}. See also \cite{DNR} (Theorem 2.2.13) or \cite{KM}.

\begin{theorem}\label{t52}
 (I)
 Assume $\vert LO(G)\vert<\infty$. Then the following holds.

(1) There is a unique rational
 series\footnote{Rational series means that for any $i$,
$G_i/G_{i+1}$ is a rational group, i.e, an abelian group embeddable into $\Q$.}
\begin{equation}\label{e51}
G=G_0\rhd G_1\rhd\cdots \rhd G_n\rhd G_{n+1}=\{e\}.
\end{equation}
(The uniqueness implies that
each subgroup $G_i$ is characteristic, i.e, invariant by any automorphism
 of $G$. Especially it is a normal subgroup of $G$.)

(2) There are elements $s_i\in G_i\setminus G_{i+1}$ for each
 $i\in\{0,1,\ldots,n\}$ such that for any map 
$\epsilon:\{0,1,\ldots,n\}\to\{\pm1\}$, there is exactly one order
 $\lambda_\epsilon$ such that $s_i^{\epsilon(i)}$ is positive. Thus
$$LO(G)=\{\lambda_\epsilon\mid \epsilon\in\{\pm1\}^{\{0,1,\ldots,n\}}\}.$$

(3) The sequence (\ref{e51}) is the maximal convex sequence for any
 $\lambda_\epsilon$. 

(4) The quotient group
 $G_i/G_{i+2}$, $i\in\{0,\ldots,n-1\}$, is not
 bi-orderable.

\medskip
(II) Conversely, if a group $G$ admits a rational series (\ref{e51})
 {  such that
 $G_{i+2}$ is a normal subgroup of $G_i$ and $G_i/G_{i+2}$ is not 
bi-orderable ($0\leq i\leq n-1$), then $\vert LO(G)\vert=2^{n+1}$.}

\end{theorem}

\bd
First of all let us show (II).
{  It suffices to prove that $G_1$ is $\lambda$-convex for any 
$\lambda\in LO(G)$. In fact, this implies that any $\lambda$ is constructed
in a lexicographical way, and thus
$\abs{LO(G)}=\abs{LO(G_1)}\cdot\abs{LO(G/G_1)}$. On the other hand, we have
$\abs{LO(G/G_1)}=2$. An induction on $n$ shows that 
$\abs{LO(G)}=2^{n+1}$.}
We use the following easy fact.

\smallskip \
\em If $A$ is a rational group and $\phi:A\to\{\pm1\}$ is a nontrivial
homomorphism, then for any nontrivial element $g\in A$, there are
$g_0\in A$ and $n\geq1$ such that $g=g_0^n$ and $\phi(g_0)=-1$.

\rm
\smallskip \rm

{  Fix $\lambda\in LO(G)$.
We shall show that $G_1$ is $\lambda$-convex by an induction on $n$.}
Consider an exact sequence
$$
1\to G_1/G_2\to G/G_2\to G/G_1\to 1.$$ 
By the induction hypothesis, $G_2$ is $\lambda\vert_{G_1}$-convex and
there is a left order $<$ on $G_1/G_2$ induced from 
$\lambda\vert_{G_1}$. 
One can define a homomorphism $\phi':G\to\{\pm1\}$
according as the conjugation by an element of $G$ preserves the order $<$ on
$G_1/G_2$ or not. (Notice that there are only two orders on $G_1/G_2$,)
Since $G_1/G_2$ is abelian, $\phi'$ induces a homomorphism
$\phi:G/G_1\to\{\pm1\}$. { Should $\phi$, equivalently $\phi'$, be
trivial, the order on $G/G_2$ constructed lexicographically from $<$
and an order of $G/G_1$ would be a bi-order. This shows that
$\phi$ is nontrivial.}

To complete the proof, let us show that for any element 
$g\in G\setminus G_1$, $g>_\lambda e$, we have 
$g^{-1}<_\lambda G_1<_\lambda g$.
There exist $g_0\in G$ and 
$n\geq 1$ such that $g\equiv g_0^n$ mod $G_1$ and $\phi'(g_0)=-1$.
Then for any $h\in G_1\setminus G_2$, $h>_\lambda e$ if and only if
$g_0^{-1}hg_0<_\lambda e$.

Assume for a while that $g_0>_\lambda e$.
Then if $h>_\lambda e$,
$$
e<_\lambda h<_\lambda hg_0<_\lambda g_0.$$
Applying $h$ successively, we obtain
$$
e<_\lambda h<_\lambda h^2<_\lambda\cdots<_\lambda h^2g_0<_\lambda hg_0<_\lambda g_0.$$
If we put $h_1=g_0^{-1}hg_0$, then
\begin{equation}\label{ex1}
e<_\lambda h<_\lambda h^2<_\lambda\cdots<_\lambda g_0h_1^2<_\lambda
g_0h_1<_\lambda g_0.
\end{equation}
By an analogous argument, we have
\begin{equation}
 \label{ex2}
g_0^{-1}<_\lambda g_0^{-1}h_1^{-1}<_\lambda
g_0^{-1}h_1^{-2}<_\lambda\cdots<_\lambda
h^{-2}<_\lambda h^{-1}<_\lambda e.
\end{equation}
The elements $h, h_1\in G_1\setminus G_2$ are
$\lambda\vert_{G_1}$-cofinal\footnote{{ For $\lambda\in
LO(G)$,
an element $h\in G$ is said to be
$\lambda$-cofinal if for any $g\in G$, there are
$n,m\in\Z$ such that $h^n<_\lambda g<_\lambda h^m$.}}
by the assumption that $G_2$ is $\lambda\vert_{G_1}$-convex and $G_1/G_2$
is rational. Therefore by (\ref{ex1}) and (\ref{ex2}), we obtain
$g_0^{-1}G_1<_\lambda G_1<_\lambda g_0G_1$.
For our initial $g$, since $g^{\pm1}G_1=g_0^{\pm n}G_1$,
we have $g^{-1}G_1<_\lambda G_1<_\lambda gG_1$, as is
required. On the other hand, if $g_0<_\lambda e$, then the same argument
shows that $gG_1<_\lambda G_1<_\lambda g^{-1}G_1$, contradicting the hypothesis
$g>_\lambda e$.
This finishes the proof of (II).
\bigskip

Now we shall proceed to the proof of (I). For a Tararin group $G$,
let $n(G)$ be the minimal height of all the elements of $LO(G)$. We shall show
(I) by the induction on $n(G)$. This is already shown for
$n(G)=0$ by Proposition 5.6. Let $G$ be a Tararin group, $\lambda\in LO(G)$
with height $n=n(G)$, and $G_1$ the maximal
proper $\lambda$-convex subgroup. Then the lexicographic 
construction shows that $G_1$ is also Tararin, and $n(G_1)\leq n-1$.
Therefore by the induction hypothesis, the maximal convex sequence of
$\lambda\vert_{G_1}$
$$
G_1\rhd G_2\rhd\cdots\rhd G_{n+1}=\{e\}$$
is a unique rational series of $G_1$ and $G_i/G_{i+2}$ is not
bi-orderable ($1\leq i\leq n-1$).

First of all, let us show that $G_1$ is a normal subgroup of $G$, and $G/G_1$
is a rational group. 
%If we show this, the elements $t_j$ can be replaced by a single element
%$s_0\in G\setminus G_1$, and the proof of the claim will be complete. 
But this is clear if the minimal set $\MM$
of the dynamical realization $\rho_\lambda$ is discrete. 
So assume $\MM$ is a locally Cantor set.
Let $x_0$ be the base point of $\rho_\lambda$, and choose $g_k\in G$ so
that $\rho_\lambda(g_k)x_0\to \exists y_0\in\MM$ as $k\to\infty$. 
One may assume that
$\rho_\lambda(g_k)x_0$ belongs to a distinct gap of $\MM$ for each $k$.
The left orders of $G$ induced by the $\rho_\lambda(G)$-orbit of $\rho_\lambda(g_k)x_0$ are
finite in number. So one may assume, by passing to a subsequence,
that the left orders are the same.
By the same argument as in Theorem
\ref{t31}, one can construct order preserving homeomorphisms $h_{k,k'}$ of
$\Cl(\rho_\lambda(G))x_0 $ which
commute with any $\rho_\lambda(g)$ such that 
$h_{k,k'}(\rho_\lambda(g_k)x_0)=\rho_\lambda(g_{k'})x_0$.
The map $h_{k,k'}$ leaves the unique minimal set $\MM$ of
$\rho_\lambda(G)$ invariant.

Consider the quotient space $\RR$ of $\R$
obtained by collapsing each gap of $\MM$ to a point.
Then $h_{k,k'}$ induces an orientation preserving homeomorphism of $\RR$
commuting with the induced action of $\rho_\lambda(G)$.
Let $Z$ be the centralizer of the action on $\RR$ induced from
$\rho_\lambda(G)$ in
the space of the orientation preserving homeomorphisms of $\RR$.
Then since  the induced action of $G$ is minimal,
$Z$ acts freely on $\RR$. In fact, if an element of $Z$ has nonempty
fixed point set, 
then the fixed point set must be $G$-invariant and coincides with $\RR$.
Thus the action of $Z$ is topologically conjugate to translations.

By the choices of $k,k'$, there are arbitrarily small translations.
That is, the action of $Z$ must be minimal.
This shows that the induced $G$-action on $\RR$ itself is also by translations.
Therefore $G_1$ is the kernel of the induced $G$-action,
  and is a normal subgroup of $G$. Finally, the left order of $G/G_1$
induced by $\lambda$ must be isolated, and hence by Theorem \ref{t31},
$G/G_1$ is rational.

\medskip

Since $G_1$ is a normal subgroup of $G$, and $G_i$ ($i\geq2$) is a
characteristic subgroup of $G_1$ by the induction hypothesis, $G_i$, especially $G_2$, is a normal
subgroup of $G$.
Next let us show that $H=G/G_2$ is not
bi-orderable. Denote $A=G_1/G_2$ and $B=G/G_1$.
There is an exact sequence
\begin{equation}\label{e713}
 1\to A \to H \stackrel{q}{\to} B\to 1.
\end{equation}
Notice that $H$ is Tararin, since otherwise lexicographic construction
would yield infinitely many left orders on $G$.
The conjugation yields a homomorphism from $H$ to ${\rm Aut}(A)$, which
projects to a homomorphism $\phi:B\to {\rm Aut}(A)$ since $A$ is abelian.
Any automorphism of $A\subset\Q$ is the multiplication by a nonzero
rational number. Thus we get $\phi:B\to\Q^\times$.
If $\phi$ takes a negative value, then $H$ does not admit a bi-order, and
we are done.
If $\phi$ is trivial, then projecting 
$H=A\times B\subset \Q^2\subset\R^2$ to $\R$ along one
dimensional linear subspaces of irrational slope yields embeddings of $H$
into $\R$, from which we obtain infinitely many left orders on $H$.
A contradiction.

Assume $\phi$ is positive valued and nontrivial.
Let $\{B_i\}$ be an exhausting increasing sequence of subgroups of $B$
which are isomorphic to $\Z$, and let $H_i=q^{-1}(B_i)$. Then the exact
sequence
$$1\to A\to H_i\to B_i\to 1$$
is split.
There is a representation $f_i:H_i\to{\rm Aff}_+(\R)$ to the group of
the orientation preserving affine transformations of the real line
such that $A$ is mapped to translations (by $A\subset\Q$ itself) and
that the split image of $B_i$ is mapped to the homotheties of ratio $\phi(B_i)$ at some
point of $\R$. Two such representations are mutually conjugate by
translations (regardless of the choice of the splittings). Therefore we
can
 arrange so that $f_{i+1}$ is an extension
of $f_i$. As the direct limit, we get a faithful representation 
$f:H\to {\rm Aff}_+(\R)$. By considering the orbit of various points of $\R$
 at which $f(H)$ acts freely, we get various left orders of $H$, 
leading to a contradiction. This finishes the proof that $H$ is
not bi-orderable.

\medskip
Finally let us show that a rational series of $G$ is unique.
By the induction hypothesis, the groups $G_i/G_{i+2}$, $1\leq i\leq n-1$,
are also not bi-orderable. So the sequence
\begin{equation}\label{seq}
G=G_0\rhd G_1\rhd\cdots \rhd G_n\rhd G_{n+1}=\{e\}
\end{equation}
satisfies the hypothesis of (II). We already know that the cardinality
of $LO(G)$ is $2^{n+1}$.

Choose $s_i\in G_i\setminus G_{i+1}$ and let $S=\{s_0,\ldots,s_n\}$.
For any $\epsilon:S\to\{\pm1\}$,
define
$$S^\epsilon=\{s_i^{\epsilon(s_i)}\mid i=0,\ldots,n\}.$$
For any $\epsilon$, we can construct a left order $\lambda^\epsilon$
 whose positive cone
contains $S^\epsilon$, lexicographically using  sequence (\ref{seq}).
Such left orders exhaust $LO(G)$, since $\vert LO(G)\vert=2^{n+1}$. 
This shows that
a rational series of $G$ is unique. In fact, any such series
gives birth to a left order lexicographically. The series is the maximal
convex sequence of that order, 
but all the $2^{n+1}$ orders have (\ref{seq}) as the maximal convex sequence. 
\qed

\begin{remark}\label{exT}
  Let $(I_i,\MM_i)$ be the $i$-th internal pair associated with the
 maximal convex sequence (\ref{e51}) of a Tararin group $G$. The next subgroup $G_{i+1}$ leaves
the gap $I_{i+1}$ of $\MM_i$ in $I_i$ invariant. But because $G_{i+1}$
 is a normal subgroup of $G_i$, it leaves all the iterates of $I_{i+1}$
under $G_i$
 invariant. By Lemma \ref{l35}, these are the only gaps of
 $\MM_i$. Therefore  $G_{i+1}$ acts trivially on $\MM_i$. That is,
there is an induced action of  $G_i/G_{i+1}$ on $\MM_i$.
If $\MM_i$ is discrete, then $G_i/G_{i+1}\cong\Z$, and the action on
 $\MM_i$ is by translation. Assume $\MM_i$ is locally Cantor.
Let $\RR_i$ be the quotient space obtained by $I_i$ by collapsing
 each gap of $\MM_i$ to a point. It is homeomorphic to $\R$.
 The quotient group $G_i/G_{i+1}$ acts on $\RR_i$ minimally and freely.
 The whole action of $G$ on $\R$ is a
 ``pileup'' of translations. Any left order is discrete if
 and only if the last group $G_n$ is isomorphic to $\Z$.
  \end{remark}

\section{Maximal convex sequence}

We shall raise one more example (other than the Tararin groups)
 of isolated orders whose height is as
big as possible.
 Let $B_n$ be the braid group of $n$ strings, with the standard
 generators $\sigma_1,\cdots,\sigma_{n-1}$. Define
$$
z_1=\sigma_1\cdots\sigma_{n-1},\ \
 z_{2}=\sigma_2\cdots\sigma_{n-1},\ \ \ldots\ ,\
z_{n-2}=\sigma_{n-2}\sigma_{n-1},\ \ z_{n-1}=\sigma_{n-1},$$
and $y_i=z_i^{(-1)^{i-1}}$. Let $P_n$ be the subsemigroup of $B_n$ generated by
$y_i$'s. Based upon a result of P. Dehornoy \cite{D}, T. V. Dubrovina and
 N. I. Dubrovin \cite{DD} have shown a remarkable fact
that $P_n\sqcup P_n^{-1}=B_n\setminus\{e\}$. The left
 order $\lambda_n$ whose positive cone is $P_n$ is called the {\em
 Dubrovina-Dubrovin order}. Since $S=\{y_1,\ldots,y_{n-1}\}$
generates $P_n$,  the order $\lambda_n$ is isolated with
characteristic positive set $S$. Moreover $\lambda_n$ can be
defined lexicographically as a twist of the Dehornoy order \cite{D},
and
 the subgroups
$$B_{n-k}^*=\langle y_{k+1},\ldots,y_{n-1}\rangle=\langle
 \sigma_{k+1},\ldots,\sigma_{n-1}\rangle$$
are $\lambda_n$-convex. 
Since $\abs{S}=n-1$, they are the only convex subgroups by Lemma \ref{l43}, and the
maximal convex sequence is given by
\begin{equation}\label{e54}
B_n>B_{n-1}^*>\cdots>B_2^*>\{e\}.
\end{equation}
 The height of $\lambda_n$ is $n-2$. The order $\lambda_n$ is discrete
 since $B^*_2\cong\Z$. The $i$-th minimal set $\MM_i$ of the $i$-th
 internal pair $(I_i,\MM_i)$ is locally Cantor, { since each term
$B_{n-k}^*$ in (\ref{e54}) is
not a normal subgroup of the previous term $B_{n-(k-1)}^*$,
because $\sigma_{k+1}\in B_{n-k}^*$, $\sigma_k\in B_{n-(k-1)}^*$,
and $\sigma_k \sigma_{k+1}\sigma_{k}^{-1}\not\in B_{n-k}^*$.}

We shall construct countably many isolated orders of $B_3$ in Section 9.
\bigskip

For an isolated order $\lambda\in LO(G)$,
we can define the {\em maximal Tararin subgroup} $G_i$ in its maximal convex
sequence
\begin{equation}
 \label{max-conv}
G>G_1>\cdots>G_n>\{e\}.
\end{equation}

For $\lambda_n$, the
maximal Tararin subgroup is $B_2^*\cong\Z$, and its height is $0$.
We shall raise questions about the isolated orders of non Tararin
groups.

\begin{question}
 Is there a non Tararin group with an isolated order whose maximal
Tararin subgroup has height $\geq1$?
\end{question}

\begin{question}
 Is there a non Tararin group with an isolated and indiscrete order?
\end{question}

There is a sufficient condition for a group to be Tararin
in terms of an isolated order on it.

\begin{proposition}
 \label{p51}
If the maximal convex sequence of an isolated order $\lambda\in LO(G)$
is subnormal,\footnote{each term is a normal subgroup of the previous term} then $G$ is a Tararin group.
\end{proposition}

\bd The proof is an induction on the height of $\lambda$. For height 0,
this is true by Proposition \ref{new}. Assume the height is $\geq1$ and consider the maximal
convex sequence of $\lambda$:
\begin{equation}\label{e53}
G=G_0\rhd G_1\rhd G_2\rhd\cdots\rhd G_n\rhd G_{n+1}=\{e\}.
\end{equation}
By the induction hypothesis, $G_1$ is a Tararin group and the
subsequence
of (\ref{e53}) that begins with $G_1$ is the unique rational series in
Theorem  \ref{t52}. Since each 
$G_i$ $2\leq i\leq n$, is a characteristic subgroup of $G_1$ and since
$G_1$ is a normal subgroup of $G$, $G_i$ is a normal subgroup of $G$.
By virtue of lemmas \ref{l41} and \ref{l42}, the order induced from
$\lambda$ on
$G/G_1$ is isolated, and of height 0. Therefore $G/G_1$ is a
rational group, by virtue of Proposition \ref{new}. That is,
the sequence (\ref{e53}) is a
rational series. 

{ 
Finally let us show that $H=G/G_2$ is not bi-orderable.
Let $A=G_1/G_2$, $B=G/G_1$ and consider the exact sequence
$$
1\to A\to H\to B\to 1.$$
As in the proof of Theorem \ref{t52}, the conjugation defines a homomorphism
$\phi:B\to{\rm Aut}(A)\subset\Q^\times$. If $\phi$ attains a negative value,
then $H$ is not bi-orderable, and we are done. 

The order $\lambda$ induces a left order $\lambda_0$ of $H$, which is the
lexicographical order given by the orders of $B$ and $A$. To
fix the idea, assume that these two orders are the natural one given by
the inclusions $B\subset\Q$ and $A\subset\Q$. Notice also that
$\lambda_0$ is isolated, since $\lambda$ is isolated.

If $\phi$ is trivial,
then $H=A\times B$. Consider the embeddings
$$
A\times B\subset\Q^2\subset\R^2.$$
Let $\pi_n:\R^2\to\R$ be the projection along an one dimensional subspace
of irrational slope $k_n$. The projection $\pi_n$ maps $A\times B$
injectively to
$\R$, and this gives a left order $\lambda_n$ of $A\times B$.
{ Then $\lambda_n\to\lambda_0$ as $k_n\downarrow0$: the $y$-coordinate
becomes more and more important as $k_n\downarrow 0$, and $\lambda_0$ is the lexicographical
order for which $B$-factor ($y$-coordinate) is of the primary importance.} Thus
$\lambda_0$ is not isolated.  

If $\phi$ is nontrivial and positive valued,
 there is an embedding $\phi$  of $H$ into ${\rm Aff}_+(\R)$
(Proof of Theorem \ref{t52}).
Points $x_n\in\R$ at which $\phi(H)$ acts freely yield left orders
$\lambda_n$ on $H$. As is observed by C. Rivas \cite{R}, we have
$\lambda_n\to\lambda_0$ as $x_n\to\infty$ (the slope
of affine transformations becomes more and more important). \qed
}

\begin{corollary}
 \label{c51}
Let $\lambda\in LO(G)$ be isolated of height 1. 
If the minimal set of the dynamical
 realization is discrete, then $G$ is a Tararin group.
\end{corollary}

\bd If the minimal set is discrete, then we get a surjective homomorphism
$\phi:G\to\Z$ and its kernel is a convex subgroup. By the previous
proposition,  $G$ is a Tararin group. \qed

\begin{example}
 The above corollary does not hold if we remove the condition that
 $\lambda$ is height 1. Let us construct an example
of isolated order $\lambda\in LO(G)$ of height 2 with discrete minimal
 set,
where $G$ is non Tararin.
 We start with the
braid group $B_3$. 
The subsemigroup $P$ generated by  $y_1=\sigma_1\sigma_2$ and
 $y_2=\sigma_2^{-1}$ is the positive cone of the Dubrovina-Dubrovin
 order $\lambda_3$. The group $B_3$ is described as
$$
B_3=\langle y_1,y_2\mid y_2y_1^2y_2=y_1\rangle.$$
There is an automorphism $\phi$ of $B_3$ which satisfies
 $\phi(y_1)=y_1^{-1}$ and $\phi(y_2)=y_2^{-1}$. Therefore if we define
a group $G$ by
$$G=\langle x, y_1,y_2\mid y_2y_1^2y_2=y_1,\  xy_1x^{-1}=y_1^{-1},\
 xy_2x^{-1}=y_2^{-1}\rangle,$$ then $B_3$ is a subgroup of $G$
 \cite{HNN}. Let $\hat P$ be the subsemigroup of $G$ generated by $x$
 and $P$. { Then we have 
$B_3=P\sqcup P^{-1}\sqcup\{e\}$, $xP=P^{-1}x$, and 
$G=\hat P\sqcup\hat P^{-1}\sqcup\{e\}$. To show the last statement, denote by
$\langle x\rangle_\pm$ the subsemigroup generated by $x^{\pm1}$.
Then $\langle x\rangle_+P^{-1}=P\langle x\rangle_+\subset \hat P$ and
 $\langle x\rangle_-P=P^{-1}\langle x\rangle_-\subset \hat P^{-1}$.
Since $B_3$ is a normal subgroup of $G$, we have
$$
G=\langle x\rangle B_3
=(\langle x\rangle_+\sqcup\langle x\rangle_-\sqcup\{e\})
(P\sqcup P^{-1}\sqcup\{e\})
$$
$$
=\langle x\rangle_+P\sqcup\langle x\rangle_-P\sqcup P
\sqcup
\langle x\rangle_+P^{-1}\sqcup\langle x\rangle_-P^{-1}\sqcup P^{-1}
\sqcup \langle x\rangle_+\sqcup \langle x\rangle_-\sqcup\{e\},$$
and each term except $\{e\}$ is contained either in $\hat P$ or in
$\hat P^{-1}$.
}

The left order $\lambda$ 
on $G$ determined by
 $\hat P$ has $B_3$ as a $\lambda$-convex normal subgroup.
In fact, 
$$B_3^{-1}x=(P\sqcup P^{-1}\sqcup\{e\})x=Px\sqcup
 P^{-1}x\sqcup\{x\}
=Px\sqcup xP\sqcup \{x\}\subset\hat P$$
and likewise $B_3^{-1}x^{-1}\subset\hat P^{-1}$, which means 
$x^{-1}<_\lambda B_3<_\lambda<x$.
 Since $G/B_3\cong\Z$,
 the minimal set associated to $\lambda$ is discrete.
The dynamics of $\lambda$ is as depicted in Figure 1.
\end{example}

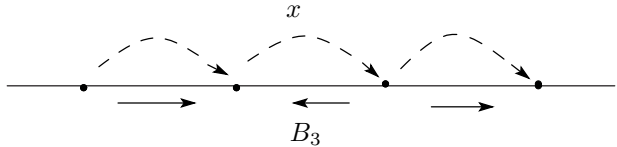
\begin{figure}[h]
%WinTpicVersion4.30a
{\unitlength 0.1in%
\begin{picture}( 31.8000,  6.7000)( 14.1000,-14.1000)%
% LINE 2 0 3 0 Black White
% 2 1410 1230 4590 1230
% 
\special{pn 8}%
\special{pa 1410 1230}%
\special{pa 4590 1230}%
\special{fp}%
% DOT 0 0 3 0 Black White
% 2 1810 1240 1810 1240
% 
\special{pn 4}%
\special{sh 1}%
\special{ar 1810 1240 16 16 0  6.28318530717959E+0000}%
\special{sh 1}%
\special{ar 1810 1240 16 16 0  6.28318530717959E+0000}%
% DOT 0 0 3 0 Black White
% 2 2610 1240 2610 1240
% 
\special{pn 4}%
\special{sh 1}%
\special{ar 2610 1240 16 16 0  6.28318530717959E+0000}%
\special{sh 1}%
\special{ar 2610 1240 16 16 0  6.28318530717959E+0000}%
% DOT 0 0 3 0 Black White
% 2 3390 1220 3390 1220
% 
\special{pn 4}%
\special{sh 1}%
\special{ar 3390 1220 16 16 0  6.28318530717959E+0000}%
\special{sh 1}%
\special{ar 3390 1220 16 16 0  6.28318530717959E+0000}%
% DOT 0 0 3 0 Black White
% 2 4190 1230 4190 1220
% 
\special{pn 4}%
\special{sh 1}%
\special{ar 4190 1230 16 16 0  6.28318530717959E+0000}%
\special{sh 1}%
\special{ar 4190 1220 16 16 0  6.28318530717959E+0000}%
% SPLINE 2 1 3 0 Black White
% 4 1890 1130 2210 980 2410 1070 2410 1070
% 
\special{pn 8}%
\special{pn 8}%
\special{pa 1890 1130}%
\special{pa 1919 1108}%
\special{pa 1941 1092}%
\special{fp}%
\special{pa 1995 1054}%
\special{pa 2006 1047}%
\special{pa 2035 1029}%
\special{pa 2050 1021}%
\special{fp}%
\special{pa 2110 994}%
\special{pa 2122 990}%
\special{pa 2151 983}%
\special{pa 2172 980}%
\special{fp}%
\special{pa 2237 985}%
\special{pa 2238 985}%
\special{pa 2267 993}%
\special{pa 2296 1005}%
\special{pa 2297 1006}%
\special{fp}%
\special{pa 2355 1036}%
\special{pa 2384 1054}%
\special{pa 2410 1070}%
\special{fp}%
% SPLINE 2 1 3 0 Black White
% 4 2670 1100 2980 970 3210 1070 3210 1070
% 
\special{pn 8}%
\special{pn 8}%
\special{pa 2670 1100}%
\special{pa 2700 1080}%
\special{pa 2725 1064}%
\special{fp}%
\special{pa 2782 1029}%
\special{pa 2788 1026}%
\special{pa 2818 1010}%
\special{pa 2841 1000}%
\special{fp}%
\special{pa 2904 978}%
\special{pa 2906 977}%
\special{pa 2936 971}%
\special{pa 2965 969}%
\special{pa 2968 969}%
\special{fp}%
\special{pa 3034 980}%
\special{pa 3053 986}%
\special{pa 3083 998}%
\special{pa 3095 1004}%
\special{fp}%
\special{pa 3154 1036}%
\special{pa 3171 1046}%
\special{pa 3200 1064}%
\special{pa 3210 1070}%
\special{fp}%
% SPLINE 2 1 3 0 Black White
% 4 3470 1140 3750 960 4010 1090 4010 1090
% 
\special{pn 8}%
\special{pn 8}%
\special{pa 3470 1140}%
\special{pa 3496 1114}%
\special{pa 3520 1091}%
\special{fp}%
\special{pa 3574 1044}%
\special{pa 3576 1042}%
\special{pa 3603 1021}%
\special{pa 3630 1002}%
\special{fp}%
\special{pa 3694 971}%
\special{pa 3711 965}%
\special{pa 3739 961}%
\special{pa 3763 961}%
\special{fp}%
\special{pa 3832 978}%
\special{pa 3852 986}%
\special{pa 3880 1001}%
\special{pa 3894 1009}%
\special{fp}%
\special{pa 3954 1048}%
\special{pa 3996 1079}%
\special{pa 4010 1090}%
\special{fp}%
% VECTOR 2 0 3 0 Black White
% 2 2460 1120 2560 1170
% 
\special{pn 8}%
\special{pa 2460 1120}%
\special{pa 2560 1170}%
\special{fp}%
\special{sh 1}%
\special{pa 2560 1170}%
\special{pa 2509 1122}%
\special{pa 2512 1146}%
\special{pa 2491 1158}%
\special{pa 2560 1170}%
\special{fp}%
% VECTOR 2 0 3 0 Black White
% 2 3270 1110 3350 1180
% 
\special{pn 8}%
\special{pa 3270 1110}%
\special{pa 3350 1180}%
\special{fp}%
\special{sh 1}%
\special{pa 3350 1180}%
\special{pa 3313 1121}%
\special{pa 3310 1145}%
\special{pa 3287 1151}%
\special{pa 3350 1180}%
\special{fp}%
% VECTOR 2 0 3 0 Black White
% 2 1990 1320 2390 1320
% 
\special{pn 8}%
\special{pa 1990 1320}%
\special{pa 2390 1320}%
\special{fp}%
\special{sh 1}%
\special{pa 2390 1320}%
\special{pa 2323 1300}%
\special{pa 2337 1320}%
\special{pa 2323 1340}%
\special{pa 2390 1320}%
\special{fp}%
% VECTOR 2 0 3 0 Black White
% 2 3200 1320 2910 1320
% 
\special{pn 8}%
\special{pa 3200 1320}%
\special{pa 2910 1320}%
\special{fp}%
\special{sh 1}%
\special{pa 2910 1320}%
\special{pa 2977 1340}%
\special{pa 2963 1320}%
\special{pa 2977 1300}%
\special{pa 2910 1320}%
\special{fp}%
% VECTOR 2 0 3 0 Black White
% 2 3630 1340 3950 1340
% 
\special{pn 8}%
\special{pa 3630 1340}%
\special{pa 3950 1340}%
\special{fp}%
\special{sh 1}%
\special{pa 3950 1340}%
\special{pa 3883 1320}%
\special{pa 3897 1340}%
\special{pa 3883 1360}%
\special{pa 3950 1340}%
\special{fp}%
% STR 2 0 3 0 Black White
% 4 2870 770 2870 870 2 0 0 0
% $x$
\put(28.7000,-8.7000){\makebox(0,0)[lb]{$x$}}%
% STR 2 0 3 0 Black White
% 4 2890 1440 2890 1540 2 0 0 0
% $B_3$
\put(28.9000,-15.4000){\makebox(0,0)[lb]{$B_3$}}%
% VECTOR 2 0 3 0 Black White
% 2 4050 1120 4140 1190
% 
\special{pn 8}%
\special{pa 4050 1120}%
\special{pa 4140 1190}%
\special{fp}%
\special{sh 1}%
\special{pa 4140 1190}%
\special{pa 4100 1133}%
\special{pa 4098 1157}%
\special{pa 4075 1165}%
\special{pa 4140 1190}%
\special{fp}%
\end{picture}}%
\caption{The dotted points form the minimal set $\MM$. The element $x$ moves
these points one to the right. The intervals bounded by the points are
invariant by $B_3$. The actions of $B_3$ are opposite in neighbouring
 intervals, showing the stability of the action.}
\end{figure}

A. Navas \cite{N} has defined the Conradian soul $C_\lambda$ for any 
$\lambda\in LO(G)$. Let us recall it briefly.
A left order $\lambda\in LO(G)$ of a group $G$ is called {\em Conradian}
if we have $g^{-1}hg^2>_{\lambda} e $ whenever
$g>_\lambda e$ and $h>_\lambda e$. Thus a bi-invariant order is
Conradian. Given an action of $G$ on $\R$, a point $x\in\R$,
is called {\em resilient} if there are an element $h$ of the stabilizer of $x$
and  a point $y\in Gx\setminus\{x\}$ such that $h^ny\to x$ as
$n\to\infty$. It is shown \cite{N}
that $\lambda\in LO(G)$ is Conradian if and only if the dynamical
realization of $\lambda$ admits no resilient point.

For a general left order $\lambda\in LO(G)$, a subgroup $H<G$ is called
$\lambda$-Conradian if the restriction of $\lambda$ to $H$ is Conradian.
The {\em Conradian soul} $C_\lambda$ of $\lambda$ is defined to be the maximal convex 
Conradian subgroup. In other words, it is the union of all the convex
Conradian subgroups.
 The following proposition is a consequence of
\cite{N}, Proposition 4.1, which states that if a group $G$ is non
Tararin,  a Conradian order of $G$ can never be isolated. Here we will give a proof based upon
Proposition \ref{p51}.

\begin{proposition}
 If $\lambda$ is isolated, the maximal Tararin subgroup of $\lambda$
 coincides with the
Conradian soul of $\lambda$. 
\end{proposition}

\bd In the maximal convex sequence (\ref{max-conv}) of $\lambda$,
let $G_i$ be the maximal Tararin subgroup. It follows from Remark
\ref{exT}, that the dynamical realization of $\lambda\vert_{G_i}$
is a pileup of translations, and cannot have a resilient point. 
Thus $G_i$ is $\lambda$-Conradian.
So it suffices to show that $G_{i-1}$ is not $\lambda$-Conradian, that
is, the dynamical realization of $\lambda\vert_{G_{i-1}}$ admits a
resilient point. It is no loss of generality to assume that $i=1$.
That is, we assume that $G$ is not a Tararin group, while its maximal convex subgroup $G_1$
is. By Proposition \ref{p51}, $G_1$ is not a normal subgroup of $G$.
Then the minimal set $\MM$ is not discrete,
and the action of $G_1$ on $\MM$ is nontrivial.
Choose $g\in G_1$ which acts nontrivially on
$\MM$. Since $G_1$ leaves invariant the gap
$I_1$ of $\MM$ containing the base point $x_0$, we have
$\Fix(g)\cap\MM\ne\emptyset$. 
Then there are distinct points $x,y\in\MM$
such that $g(x)=x$ and either $g^n(y)\to x$ or $g^{-n}(y)\to x$ as $n\to\infty$.
Since the action of $G$ on $\MM$ is minimal, the point $y$ is
accumulated by the orbit of $x$. This shows that the point $x$ is
resilient. \qed

\section{Circular orders}

In this section, we provide preliminary facts about circular orders.

\begin{definition}\label{co}
 For a countable group $\overline{G}$, a map $c:\overline{G}^3\to\{0,1,-1\}$ is called a
{\em left invariant circular order of $\overline{G}$} if it satisfies the following 
conditions.

(1) $c(g_1,g_2,g_3)=0$ if and only if $g_i=g_j$ for some $i\ne j$.

(2) For any $g_1,g_2,g_3,g_4\in \overline{G}$, we have
$$c(g_2,g_3,g_4)-c(g_1,g_3,g_4)+c(g_1,g_2,g_4)-c(g_1,g_2,g_3)=0.$$

(3) For any $g_1,g_2,g_3, g_4\in \overline{G}$, we have
$$
c(g_4g_1,g_4g_2,g_4g_3)=c(g_1,g_2,g_3).$$
\end{definition}

\begin{definition}
 Given a finite 
set $F$ of $\overline{G}$, a {\em configuration of $F$ in $S^1$} is
an equivalence class of injections $\iota:F\to S^1$, where two
 injections $\iota$ and $\iota'$ is said to be equivalent if there is
an orientation preserving homeomorphism $h$ of $S^1$ such that $\iota'=h\iota$.
\end{definition}

Given a left invariant circular order $c$ of $\overline{G}$, the
configuration of the set $\{g_1,g_2,g_3\}$ of three points
is determined by the rule that $g_1,g_2,g_3$ is positioned
anticlockwise if $c(g_1,g_2,g_3)=1$, and clockwise if $c(g_1,g_2,g_3)=-1$.
By condition (2) of Definition \ref{co}, this is well defined.
But (2) says more. One can show the following proposition by an easy
induction on the cardinality of $F$. 

\begin{proposition}\label{should}
 Given a left invariant circular order of $\overline{G}$, the configuration of any
 finite set $F$ in $S^1$ is determined. \qed
\end{proposition}

Denote by $CO(\overline{G})$ the set of all the left invariant circular orders.
It is equipped with a totally disconnected compact metrizable topology,
just as $LO(\overline{G})$. An {\em isolated left invariant circular
order} is defined using this topology.
If $c\in CO(\overline G)$ is isolated, then there is a finite set
$\overline S$ of $\overline G$, called a {\em determining set}, such that any left
invariant circular order  which gives the same configuration of $\overline S$
as $c$ is $c$.

Given $c\in CO(\overline{G})$, we define a dynamical realization
$\rho_c:\overline{G}\to {\rm Homeo}_+(S^1)$ based at $y_0\in S^1$ as follows.
Fix an enumeration of $\overline{G}$: $\overline{G}=\{g_i\mid i\in\N\}$ such that $g_1=e$.
Define an embedding $\iota:\overline{G}\to S^1$ inductively as follows. 
First, set $\iota(g_1)=y_0$ and $\iota(g_2)=y_0+1/2$. If $\iota$ is
defined on $\{g_1,\cdots,g_n\}$, then there is a connected component of 
$S^1\setminus\{\iota(g_1),\ldots,\iota(g_n)\}$ where the point $g_{n+1}$
should be embedded, by virtue of Proposition \ref{should}. Define
$\iota(g_{n+1})$ to be the midpoint of that interval.
Using the injection $\iota$, we can define the action of $\overline{G}$
on $S^1$
just as in the case of left orders. 
The action is called the {\em dynamical realization of $c$ based at $y_0$} and denoted by
$\rho_c$. We shall raise fundamental
properties of $\rho_c$. The proof is completely
parallel to the case of left orders.

\begin{lemma}\label{tight}
 The dynamical realization $\rho_c$ is {\em tight} at the base point
$y_0$, i.e, it is free at $y_0$ and if $I$ is a connected component of
$S^1\setminus{\rm Cl}(\rho_c(\overline G)y_0)$, then 
$\partial I\subset \rho_c(\overline G)y_0$.
\end{lemma}

\begin{lemma}\label{conj}
 Two dynamical realizations obtained via different enumerations of $\oG$
 are mutually conjugate by an orientation and base point preserving
 homeomorphism of $S^1$.
\end{lemma}

Let $\MM$ be a minimal set of the dynamical realization $\rho_c$ of an
isolated circular order $c$.
It is shown by K. Mann and C. Rivas \cite{MR} that (unlike left orders)
 $\MM$ is always a proper subset of $S^1$. Summarizing with other
 properties, we get:

\begin{lemma}\label{p71}
 If $\overline G$ is not finite cyclic, the minimal set $\MM$ of the dynamical
 realization $\rho_c$ of any isolated circular order $c\in CO(\oG)$ is unique.
It is either a finite set or a Cantor set.
\end{lemma}

{  \bd If $\MM$ is a Cantor set, then a standard argument shows that 
it is the unique minimal set. So suppose $\MM$ is a finite set of
cardinality, say $n$. Then by looking at the action of $\oG$ on $\MM$,
one can define a surjective homomorphism $\phi: \oG\to \Z/n\Z$.
Since $\oG$ is not finite cyclic, the base point $y_0$ of the dynamical
realization $\rho_c$ must be contained in a gap, say $I_1=(a,b)$
\footnote{Given two points $a,b\in S^1$, we define
$(a,b)=\{t\in S^1\mid a\prec t\prec b\}$, where $\prec$ is the
anticlockwise circular order of $S^1$.} of $\MM$.
The stabilizer $\oG_1$ of $I_1$ is the kernel of $\phi$, and 
$\rho_c(\oG)y_0\cap I_1=\rho(\oG_1)y_0$.
Now if we set $x=\inf(\rho_c(\oG_1)y_0)$, then $x=a$: for,
otherwise, $x\in \rho_c(\oG)y_0$ by the tightness, contradicting that $x$ is
the infimum. Likewise we have
$\sup(\rho_c(\oG_1)y_0)=b$.
This shows that all the orbits other than $\MM$ is infinite and contains
$\MM$ in its closure, that is, 
$\MM$ is the unique minimal set.
\qed

\bigskip
Taking into account the Cantor minimal set case, the similar argument
shows
the following.}

\begin{lemma}\label{p72}
 If $\oG$ is not finite cyclic and $c$ is isolated, 
then the base point $y_0$ of the dynamical
 realization is contained in a gap $I$ of the minimal set $\MM$, the
stabilizer $\oG_{I}$ of $I$ is nontrivial, and there is no gap of $\MM$
 other than the orbit of $I$.
\end{lemma}

{ Here is an analogue of Corollary 2.4 for circular orders.

\begin{lemma}
 \label{2.4.c}
Let $\overline\HH$ be the set of the orientation and base point $y_0$
 preserving topological conjugacy classes of the homomorphisms
 $\overline G\to{\rm Homeo}_+(S^1)$ which are tight at $y_0$. Then the
dynamical realization at $y_0$ induces a bijection of $CO(\overline G)$ onto
$\overline\HH$.
\end{lemma}
}

{ 
\begin{definition}
 Let $c$ be a circular order of $\oG$, isolated or not, and $H$ a
 nontrivial subgroup of $\oG$. 
$H$ is said to be {\em $c$-convex} if $\rho_c(H)$ acts
with global fixed points, and 
$\rho_c(\oG)y_0\cap I_H=\rho_c(H)y_0$, 
where $I_H$ denotes the connected component of the complement of the
global fixed point set of $\rho_c(H)$ containing $y_0$.
The configuration of $\rho_c(H)y_0$ in $I_H$ defines a left order
$\lambda$ on $H$, which we call the left order on $H$ {\em induced from
 $c$}. The trivial subgroup is said to be $c$-convex.
\end{definition}

As shown in \cite{MR} Lemma 3.15, there is a unique maximal $c$-convex
subgroup
for any $c\in CO(\oG)$, which we call the {\em linear part of $c$.}
By virtue of Lemmas \ref{p71} and \ref{p72}, we get the following lemma.

\begin{lemma}\label{linear-part}
Assume $\oG$ is not finite cyclic. Let $\MM$ be the minimal set of the
 dynamical realization of an isolated circular order $c\in CO(\oG)$, and $I$ the gap of $\MM$ which contains
 the base point $y_0$. Then the linear part of $c$ coincides with the
 stabilizer
of $I$. 
\end{lemma}
}

\section{Isolated left orders on $B_3$}

In this section, using a method of \cite{MR},
 we construct countably many isolated left orders on the braid
group $B_3$,  which
are not the automorphic images of the others. 
 The group $B_3$ has the following presentations.
\begin{equation*}
\begin{array}{lll}
 B_3&=&\langle \sigma_1,\sigma_2\mid
  \sigma_1\sigma_2\sigma_1=\sigma_2\sigma_1\sigma_2 \rangle\\
 &=&\langle y_1,y_2\mid y_2y_1^2y_2=y_1\rangle\\
&= &\langle a, b\mid a^2=b^3\rangle,
\end{array}
\end{equation*}
where the generators are related by
{ 
$$
y_1=\sigma_1\sigma_2,\ \ y_2=\sigma_2^{-1},\ \ a=y_2y_1^2,\ \ b= y_1.$$
The Dubrovina-Dubrovin order
$\lambda_3$ 
is the unique left order on $B_3$
which satisfies $y_1>_{\lambda_3}e$ and $y_2>_{\lambda_3}e$,
equivalently 
$e<_{\lambda_3}a<_{\lambda_3}b$.
To show the equivalence, assume $y_1>_{\lambda_3}e$ and $y_2>_{\lambda_3}e$. 
Then 
$$a=y_2y_1^2>_{\lambda_3}e,\ \mbox{ and }\
a^{-1}b=y_1^{-2}y_2^{-1}y_1=y_1^{-2}y_2^{-1}(y_2y_1^2y_2)=y_2>_{\lambda_3}e.$$
The converse is shown similarly.}

Henceforth in this section we denote by $G$ the braid group $B_3$ and 
by $\oG$ its
quotient by the center. Namely, we put
$$
G=\langle a,b,t\mid a^2=b^3=t\rangle,\ \ \
\overline G=\langle\alpha,\beta\mid\alpha^2=\beta^3=e\rangle,$$
and  $q:G\to\overline G$ to be the surjective homomorphism satisfying
$q(a)=\alpha$ and $q(b)=\beta$. 

{  The first half of this section is devoted to show that
there is a homeomorphism between $LO(G)$ and $CO(\oG)$. Thus the
construction of isolated orders in $LO(G)$ reduces to the construction
of isolated orders in $CO(\oG)$, which is easier, thanks to well
developed theory of Fuchsian groups. The last half is devoted to this
construction.

First of all, notice that $t$ is $\lambda$-cofinal for any 
$\lambda\in LO(G)$. In fact, assume, to fix the idea, that $t>_\lambda
e$.
Then $a>_\lambda e$ and $b>_\lambda e$. Since $t$ is in the center of $G$,
any element $g\in G$ can be written as
$g=t^Na^{i_1}b^{j_1}\cdots a^{i_r}b^{j_r}$ for an integer $N$ and
 negative integers
$i_\nu$ and $j_\nu$, showing that $g<_\lambda t^N$. Likewise
there is $M\in\Z$ such that $t^M<_\lambda g$. }

Let $\tau$ be the translation
of $\R$ by 1, and  $p:\R\to S^1=\R/\langle \tau\rangle$ the canonical
projection. Denote by ${\rm Homeo}_\Z(\R)$ the group of all the
homeomorphisms of $\R$ which commute with $\tau$.
It is the universal covering group of ${\rm Homeo}_+(S^1)$.
Denote by
$\pi:{\rm Homeo}_\Z(\R)\to{\rm Homeo}_+(S^1)$ the covering map.
{ Let 
$$LO^+(G)=\{\lambda\in LO(G)\mid t>_\lambda e\}\ \mbox{ and }\
LO^-(G)=\{\lambda\in LO(G)\mid t<_\lambda e\}.$$
Likewise let $$CO^\pm(\oG)=\{c\in CO(\oG)\mid c(e,\beta,\beta^2)=\pm1\}.$$
In order to show that $LO^\pm(G)$ is homeomorphic to $CO^\pm(\oG)$, it is
easier and more natural to consider the sets of the conjugacy classes of
certain homeomorphisms: one is $\HH^\pm$ with a bijection to $LO^\pm(G)$ and
the other is
$\overline\HH^\pm$ with a bijection to $CO^\pm(\oG)$. Next we shall
 construct a natural
bijection from $\HH^\pm$
to $\oH^\pm$. Finally, we shall show that the induced
bijection from $LO^\pm(G)$ to $CO^\pm(\oG)$ is a homeomorphism
(Theorem 9.1).

Let
 $H^\pm$ be the set of the 
homomorphisms $\rho:G\to {\rm Homeo}_\Z(\R)$ which are tight at
a prescribed base point
$x_0$ and satisfy $\rho(t)=\tau^{\pm1}$, and let $\HH^\pm$ be the
set of the orientation and base point preserving topological conjugacy classes
of the elements of $H^\pm$.
{  Then we have $\HH=\HH^+\cup\HH^-$ for $\HH$ in Corollary
 2.4. In fact, for any class $[\rho]$ of $\HH$, $\rho(t)$ is fixed point
 free since $t$ is cofinal for any left orders. Thus $[\rho]$ has a
representative $\rho$ such that $\rho(t)=\tau^{\pm1}$.
}

A dynamical realization at $x_0$ of any element of $LO^{\pm}(G)$
represents a unique element of $\HH^\pm$
(Corollary \ref{bijection}).
That is, we get a map $\phi:LO^{\pm}(G)\to\HH^{\pm}$. 
Clearly $\phi$ is injective. To show the surjectivity, let
$\rho\in H^\pm$, let $\lambda\in LO^\pm(G)$ be the left order of $G$
defined by the natural order of the orbit $\rho(G)x_0$,
and let $\rho_\lambda\in H^\pm$ a representative of the conjugacy
class of the dynamical realization of $\lambda$. There is an order
and the base point
preserving equivariant bijection from $\rho_\lambda(G)x_0$  to
$\rho(G)x_0$, which can be extended to a homeomorphism between the
closures, thanks to the tightness. Finally this homeomorphism can be extended to an equivariant
homeomorphism of $\R$. We have completed the proof that $\phi$ is a
bijection.

Denote by $\overline H^\pm$ the set of the 
homomorphisms 
$\overline\rho:\overline G\to{\rm Homeo}_+(S^1)$ which are tight at
$y_0=p(x_0)$ and satisfy ${\rm rot}(\overline\rho(\beta))=\pm1/3$, and let
$\overline\HH^{\pm}$ the set of the
orientation and the base point preserving
 topological conjugacy classes  of the  elements of $\overline H^\pm$.
 Then 
$\overline \HH^\pm$ is identified with
$CO^\pm(\oG)$ (Lemma \ref{2.4.c}).

}
Define a map $q_*:H^{\pm}\to\overline H^{\pm}$ by 
$(q_*\rho)(\overline g)=\pi(\rho(g))$, where $\rho\in H^\pm$, 
$\overline g\in\oG$ and $g\in G$ is any element such that 
$q(g)=\overline g$. 
There is a commutative diagram
$$\begin{array}{ccl}
G & \stackrel{\rho}{\to} & {\rm Homeo}_\Z(\R) \\
\downarrow q & & \ \ \ \ \downarrow\pi\\
\oG & \stackrel{q_*\rho}{\to} & {\rm Homeo}_+(S^1).
\end{array}$$
Define a map $\pi^*:\overline H^\pm\to H^\pm$ for
$\overline\rho\in \overline H^\pm$ by 

$\bullet$ $(\pi^*\overline\rho)(a)$ is the lift of
	    $\overline\rho(\alpha)$ to ${\rm Homeo}_\Z(\R)$ whose square
	    is
            $\tau^{\pm1}$, and

 $\bullet$ $(\pi^*\overline\rho)(b)$ is the lift of
	    $\overline\rho(\beta)$ to ${\rm Homeo}_\Z(\R)$ whose cube
	    is
            $\tau^{\pm1}$.
\\ 
Also we have a commutative diagram
$$\begin{array}{ccl}
G & \stackrel{\pi^*\overline\rho}{\to} & {\rm Homeo}_\Z(\R) \\
\downarrow q & & \ \ \ \ \downarrow\pi\\
\oG & \stackrel{\overline\rho}{\to} & {\rm Homeo}_+(S^1).
\end{array}$$
It is clear that $q_*$ and $\pi^*$ map the conjugacy classes to the
conjugacy classes. That is, we have maps (denoted by the same letters)
$q_*:\HH^\pm\to\oH^\pm$ and $\pi^*:\oH^\pm\to\HH^\pm$.
We have $\pi^*q_*=q_*\pi^*=id$. Thus we obtain a bijection 
$q_*:LO(G)\to CO(\oG)$ and its inverse $\pi^*$.

We shall show the following theorem (Theorem \ref{t4} in the introduction).

\begin{theorem}\label{t91}
 The map $q_*:LO(G)\to CO(\oG)$ is a homeomorphism.
\end{theorem}

\bd
 Let us show that $q_*$
is continuous.
For any $\lambda\in LO(G)$,
 let $c=q_*(\lambda)\in CO(\oG)$. Choose arbitrary elements $\og_1,\ldots,\og_n$ of $\oG$ and
consider their configuration in $S^1$ with respect to $c$. This is the same as the configuration
of $\overline\rho(\og_1)y_0,\ldots,\overline\rho(\og_n)y_0$ in $S^1$, where 
$\overline\rho\in\overline\HH$ is a dynamical realization of $c$. Let
$\rho=\pi^*(\overline\rho)\in\HH$, a dynamical realization of $\lambda$.
Choose  $g_i\in G$ 
such that $q(g_i)=\og_i$ and $e\leq_\lambda g_i<_\lambda t$ ($1\leq i\leq n$).
The configuration of 
$\rho(g_1)x_0,\ldots,\rho(g_n)x_0$ in $\R$
coincides with the
configuration of $g_1,\cdots,g_n$ with respect to $\lambda$. Choose any
$\lambda'\in LO(G)$ whose configuration of $e,g_1,\cdots,g_n,t$ is the
same as $\lambda$. Then the configuration of $\og_1,\cdots,\og_n$ of 
$q_*(\lambda')$ is the same as $c$, showing the continuity of $q_*$.
Thus the compact metrizable sets $LO(G)$ and $CO(\oG)$ are homeomorphic by $q_*$.

\qed

By virtue of the previous theorem, Theorem \ref{t3} in the introduction
reduces to the following theorem. This is because any automorphism of
$G$,
{  preserving the center,}
induces an automorphism of $\oG$.

\begin{theorem}\label{t92}
 There are isolated circular orders $c^{(k)}\in CO(\oG)$, ($k>0,k\equiv\pm1$mod 6) which
 are not the automorphic images of the others.
\end{theorem}

The rest of this section is devoted to the proof of this theorem.
{ Our argument is based upon Fuchsian representations.
The Lie group $PSL(2,\R)$ is the group of the orientation preserving
isometries of the Poincar\'e upper half plane 
$\mathbb H=\{z\in\mathbb C\mid {\rm Im}(z)>0\}$, acting by linear fractional
transformations. We consider $\mathbb H$ to be an open half disk
in the Riemann sphere $\hat{\mathbb C}=\mathbb C\cup\{\infty\}$.
Then $PSL(2,\R)$ acts on $\hat{\mathbb C}$ by linear fractional transformations
and leaves invariant the oriented boundary $\partial\mathbb H$, which we
identify with $S^1$. This gives an inclusion 
$PSL(2,\R)\subset{\rm Homeo}_+(S^1)$. Denote by
$\overline{\mathbb H}$ the closure of $\mathbb H$ in $\hat{\mathbb C}$.

Let $\Gamma$ be a nonamenable countable group. A representation
$\rho:\Gamma\to PSL(2,\R)$ is called {\em Fuchsian} if it is faithful
and the image $\rho(\Gamma)$ is discrete in $PSL(2,\R)$. For a Fuchsian
representation $\rho$, its {\em limit set} 
$L_\rho\subset \partial \mathbb H\approx S^1$ is, by definition, the set
of the accumulation points of an
orbit $\rho(\Gamma)z_0$
($z_0\in\mathbb H$). It does not depend on the choice of
$z_0$. It is also characterized as the unique minimal set of the
representation 
$\Gamma\stackrel{\rho}{\to}PSL(2,\R)\subset {\rm Homeo}_+(S^1)$.
It is either the whole $S^1$ or a Cantor set (by the nonamenability of 
$\Gamma$). In the former case, $\rho$ is called  {\em of the first
kind}, and in the latter {\em of the second kind}.
 
}

Let us return to our group $\oG$.
There is an isomorphism $\iota:\overline G \cong PSL(2,\Z)$ which satisfies
$$
\iota(\alpha)=\left[\begin{array}{cc}0&-1\\1&0\end{array}\right],\ \ 
\iota(\beta)=\left[\begin{array}{cc}1&1\\ -1&0\end{array}\right].$$

Let us define a homomorphism $\rho_M:\overline G\to PSL(2,\R)$, called the {\em modular representation},
as the composite
$$
\overline G\stackrel{\iota}{\cong}PSL(2,\Z)\subset PSL(2,\R).
$$ 
For the dynamics of the modular representation
$\oM$, see Figure 2. The open disk bounded by the circle is
{  the image of $\mathbb H$ by the stereographic projection
from $-i$.}
\begin{figure}[h]
%WinTpicVersion4.30a
{\unitlength 0.1in%
\begin{picture}( 22.4000, 20.9900)( 15.0000,-22.3000)%
% DOT 0 0 3 0 Black White
% 2 2635 1269 2635 1269
% 
\special{pn 4}%
\special{sh 1}%
\special{ar 2635 1269 16 16 0  6.28318530717959E+0000}%
\special{sh 1}%
\special{ar 2635 1269 16 16 0  6.28318530717959E+0000}%
% CIRCLE 2 0 3 0 Black White
% 4 2640 1269 3394 696 3394 696 3296 784
% 
\special{pn 8}%
\special{ar 2640 1269 947 947  5.6465485  5.6333485}%
% LINE 2 0 3 0 Black White
% 2 2640 328 2640 2220
% 
\special{pn 8}%
\special{pa 2640 328}%
\special{pa 2640 2220}%
\special{fp}%
% DOT 0 0 3 0 Black White
% 2 2145 1269 2145 1269
% 
\special{pn 4}%
\special{sh 1}%
\special{ar 2145 1269 16 16 0  6.28318530717959E+0000}%
\special{sh 1}%
\special{ar 2145 1269 16 16 0  6.28318530717959E+0000}%
% SPLINE 2 0 3 0 Black White
% 4 2640 323 2424 897 2150 1269 2150 1269
% 
\special{pn 8}%
\special{pa 2640 323}%
\special{pa 2613 416}%
\special{pa 2603 447}%
\special{pa 2594 478}%
\special{pa 2574 540}%
\special{pa 2564 570}%
\special{pa 2554 601}%
\special{pa 2521 691}%
\special{pa 2509 721}%
\special{pa 2497 750}%
\special{pa 2484 779}%
\special{pa 2470 808}%
\special{pa 2456 836}%
\special{pa 2442 865}%
\special{pa 2427 892}%
\special{pa 2411 920}%
\special{pa 2377 974}%
\special{pa 2341 1026}%
\special{pa 2284 1104}%
\special{pa 2204 1204}%
\special{pa 2162 1254}%
\special{pa 2150 1269}%
\special{fp}%
% SPLINE 2 0 3 0 Black White
% 4 2140 1259 2453 1622 2640 2220 2640 2220
% 
\special{pn 8}%
\special{pa 2140 1259}%
\special{pa 2163 1282}%
\special{pa 2187 1305}%
\special{pa 2233 1351}%
\special{pa 2255 1374}%
\special{pa 2278 1398}%
\special{pa 2300 1421}%
\special{pa 2342 1469}%
\special{pa 2363 1494}%
\special{pa 2401 1544}%
\special{pa 2420 1570}%
\special{pa 2437 1596}%
\special{pa 2469 1650}%
\special{pa 2497 1706}%
\special{pa 2510 1735}%
\special{pa 2522 1764}%
\special{pa 2533 1793}%
\special{pa 2544 1823}%
\special{pa 2554 1854}%
\special{pa 2563 1884}%
\special{pa 2581 1946}%
\special{pa 2588 1978}%
\special{pa 2596 2009}%
\special{pa 2624 2137}%
\special{pa 2630 2170}%
\special{pa 2636 2202}%
\special{pa 2640 2220}%
\special{fp}%
% LINE 2 0 3 0 Black White
% 2 2150 1269 1689 1269
% 
\special{pn 8}%
\special{pa 2150 1269}%
\special{pa 1689 1269}%
\special{fp}%
% CIRCLE 3 0 3 0 Black White
% 4 2635 1269 2752 1279 2679 1387 2576 1161
% 
\special{pn 4}%
\special{ar 2635 1269 117 117  4.2123938  1.2138844}%
% VECTOR 3 0 3 0 Black White
% 2 2586 1171 2542 1191
% 
\special{pn 4}%
\special{pa 2586 1171}%
\special{pa 2542 1191}%
\special{fp}%
\special{sh 1}%
\special{pa 2542 1191}%
\special{pa 2611 1182}%
\special{pa 2591 1169}%
\special{pa 2594 1145}%
\special{pa 2542 1191}%
\special{fp}%
% CIRCLE 3 0 3 0 Black White
% 4 2150 1264 2311 1298 2208 1097 1988 1210
% 
\special{pn 4}%
\special{ar 2150 1264 165 165  3.4633432  5.0466612}%
% VECTOR 3 0 3 0 Black White
% 2 1993 1200 1983 1230
% 
\special{pn 4}%
\special{pa 1993 1200}%
\special{pa 1983 1230}%
\special{fp}%
\special{sh 1}%
\special{pa 1983 1230}%
\special{pa 2023 1173}%
\special{pa 2000 1179}%
\special{pa 1985 1160}%
\special{pa 1983 1230}%
\special{fp}%
% CIRCLE 2 0 3 0 Black White
% 4 2648 1269 3691 961 3705 1003 3299 408
% 
\special{pn 8}%
\special{ar 2648 1269 1088 1088  5.3597811  6.0366490}%
% STR 2 0 3 0 Black White
% 4 2592 191 2592 261 2 0 0 0
% $\infty$
\put(25.9200,-2.6100){\makebox(0,0)[lb]{$\infty$}}%
% STR 2 0 3 0 Black White
% 4 1500 1248 1500 1318 2 0 0 0
% $-1$
\put(15.0000,-13.1800){\makebox(0,0)[lb]{$-1$}}%
% STR 2 0 3 0 Black White
% 4 2081 1304 2081 1374 2 0 0 0
% $\omega$
\put(20.8100,-13.7400){\makebox(0,0)[lb]{$\omega$}}%
% STR 2 0 3 0 Black White
% 4 2515 1262 2515 1332 2 0 0 0
% $i$
\put(25.1500,-13.3200){\makebox(0,0)[lb]{$i$}}%
% VECTOR 2 0 3 0 Black White
% 2 3710 990 3730 1070
% 
\special{pn 8}%
\special{pa 3710 990}%
\special{pa 3730 1070}%
\special{fp}%
\special{sh 1}%
\special{pa 3730 1070}%
\special{pa 3733 1000}%
\special{pa 3717 1018}%
\special{pa 3694 1010}%
\special{pa 3730 1070}%
\special{fp}%
% STR 2 0 3 0 Black White
% 4 3730 460 3730 560 2 0 0 0
% $\rho_M(\alpha\beta)$
\put(37.3000,-5.6000){\makebox(0,0)[lb]{$\rho_M(\alpha\beta)$}}%
% STR 2 0 3 0 Black White
% 4 2860 1230 2860 1330 2 0 0 0
% $\rho_M(\alpha)$
\put(28.6000,-13.3000){\makebox(0,0)[lb]{$\rho_M(\alpha)$}}%
% STR 2 0 3 0 Black White
% 4 1870 980 1870 1080 2 0 0 0
% $\rho_M(\beta)$
\put(18.7000,-10.8000){\makebox(0,0)[lb]{$\rho_M(\beta)$}}%
% STR 2 0 3 0 Black White
% 4 2590 2260 2590 2360 2 0 0 0
% $0$
\put(25.9000,-23.6000){\makebox(0,0)[lb]{$0$}}%
\end{picture}}%

\caption{}
\end{figure}
The element $\oM(\alpha)$ is the $1/2$-rotation
around $i$, and $\oM(\beta)$ the $1/3$-rotation around
$\omega=(-1+\sqrt{-3})/2$. The element $\oM(\alpha\beta)$ is a parabolic
transformation which fixes the point $0$ and moves points on
$S^1\setminus\{0\}$ clockwise, as is depicted in Figure 2.
The Fuchsian group $\oM(\oG)$ is of the first kind.

Let us define another Fuchsian representation
$\rho:\overline G\to PSL(2,\R)$, a deformation of $\oM$. Choose a point
$\omega'$ on the geodesic which passes through $i$ and $\omega$, but
slightly farther than $\omega$ from $i$: $d(\omega',i)>d(\omega,i)$.
See Figure 3. We set $\rho(\alpha)$ to be the same as $\oM(\alpha)$, the
$1/2$-rotation around $i$, and $\rho(\beta)$ the $1/3$-rotation around
$\omega'$. We put the base point $y_0=0$.

{ Consider the 4-gon $P$ depicted in Figure 3, a closed subset of 
$\overline{\mathbb H}$. It is routine to show, using hyperbolic metric
of $\mathbb H$, that the translates of $P\cap\mathbb H$ tesselate 
$\mathbb H$, that is, $\bigcup_{\gamma\in\oG}\rho(\gamma)(P\cap\mathbb H)=\mathbb H$ and
$\rho(\gamma)({\rm Int}P)\cap{\rm Int}P=\emptyset$ if $\gamma\neq e$.
This shows that $\rho$ is a Fuchsian representation. Moreover it is of
the second kind since the tesselation implies that $L_\rho\cap P=\emptyset$.
But we need a bit more: we shall show that the translates of $P$ 
tesselate the whole $\overline{\mathbb H}\setminus L_\rho$, that is,
$\bigcup_{\gamma\in\oG}\rho(\gamma)P=\overline{\mathbb H}\setminus L_\rho$.

Endow the half disk $\overline{\mathbb H}$ with the restriction
of the spherical metric of the Riemann sphere 
$\mathbb C\cup\{\infty\}$. Then for any $g\in PSL(2,\R)$ and any
$z\in\overline{\mathbb H}$, the absolute value of the derivative
$\abs{g'(z)}$ is well defined. Given $g\in PSL(2,\R)\setminus PSO(2)$, define
$I(g)=\{z\in\overline{\mathbb H}\mid\abs{g'(z)}\geq 1\}$.
It is a subset of $\overline{\mathbb H}$ delimited by a circle
perpendicular to $\partial{\mathbb H}$, whose radius is
denoted by ${\rm rad}(I(g))$.
It satisfies the following properties: (1) 
$gI(g)=\overline{\overline{\mathbb H}\setminus I(g^{-1})}$.
(2) $g$ is hyperbolic if and only if\footnote{To show ``only if part'', let $\ell$
be the axis of a hyperbolic element $g$. Let $m$ be the midpoint of
$\ell$
and let $p,q\in\ell$ be points such that $d(p,m)=d(q,m)$ and
$g(p)=q$. Then $I(g)$ (resp.\ $I(g^{-1})$) is bounded by a circle crossing $\ell$
perpendicularly at $p$ (resp.\ $q$). Thus $I(g)\cap I(g^{-1})=\emptyset$.}  $I(g)\cap I(g^{-1})=\emptyset$.
(3) $g_n\to \infty$ in $PSL(2,\R)$ if and only if ${\rm rad}(I(g_n))\to 0$.

Returning to our representation $\rho:\oG\to PSL(2,\R)$, there is a
purely hyperbolic subgroup of $\rho(\oG)$ of index 6, as is shown 
later. On the other hand, the limit set does not change if we pass
to a finite index subgroup. Thus we get:

\smallskip \em
For a point $x\in\partial\mathbb H$, $x\in L_\rho$ if and only if
$\rho(\gamma_n)z_0\to x $ for some hyperbolic elements
$\rho(\gamma_n)$, where $z_0$ is any prescribed point in $\mathbb H$.

\smallskip\rm
Now if $\rho(\gamma)$ is hyperbolic, then $I(\rho(\gamma))$, as well as
$I(\rho(\gamma)^{-1})$, contains a fixed point of $\rho(\gamma)$. 
In particular,
$I(\rho(\gamma))\cap L_\rho\neq\emptyset$ and
$I(\rho(\gamma)^{-1})\cap L_\rho\neq\emptyset$.
Moreover if $\rho(\gamma)$ is sufficiently far away from the identity,
then $I(\rho(\gamma))\cap P=\emptyset$ because 
${\rm rad}(I(\rho(\gamma)))$ is small and $P\cap L_\rho=\emptyset$.
Therefore the translate $\rho(\gamma)P$ is contained in
$I(\rho(\gamma)^{-1})$ and has small diameter. This is true not only for
elements of a subgroup of index 6, but also for any element of
$\rho(\oG)$: if $\gamma_n\to\infty$ in the word norm, then
${\rm diam}(\rho(\gamma_n)P)\to 0$.
This yields yet another characterization of the limit set.

\smallskip\em
For $x\in\partial\mathbb H$, $x\in L_\rho$ if and only if
any neighbourhood of $x$ in $\overline{\mathbb H}$ intersects
infinitely many translates of $P$.

\smallskip\rm

Finally, this, together with the fact that the translates of
$P\cap\mathbb H$ tesselate $\mathbb H$, implies that the translates of $P$ tesselate
$\overline{\mathbb H}\setminus L_\rho$. This is what we wanted to show.

\medskip
Henceforth we denote the composite
$$
\oG \stackrel{\rho}{\to}PSL(2,\R)\subset{\rm Homeo_+}(S^1)$$ 
also by $\rho$,
and  the unique minimal set $L_\rho$ of the
homomorphism $\rho$ by $\MM$.
The translates of the interval $P\vert_{S^1}$ tesselates
$S^1\setminus\MM$,
where $P\vert_{S^1}$ is the component of $P\cap S^1$ which is
homeomorphic to a closed interval.
Let $y_0$ be a point in $S^1$ which corresponds to $0\in\partial\mathbb H$
and is depicted by $e$ in Figure 3. Since $y_0\in P$, it lies in a gap
of $\MM$, say $I_1=(\sigma_-,\sigma_+)$. The element $\alpha\beta\in\oG$ 
sends $y_0$ to a point slightly right to itself ($e$ to $\alpha\beta$ in
Figure 3). 
As can be seen by Figure 3,
{  the iterates of 
$\rho(\alpha) P\vert_{S^1}$ by $\rho(\alpha\beta)^n$, $n\in\Z$, tesselate an open
interval
of $S^1\setminus\MM$. Therefore
$\rho(\alpha\beta)^ny_0$ lie in $I_1$.} We also have 
$\lim_{n\to\pm\infty}\rho(\alpha\beta)^ny_0=\sigma_\pm$, since the
limits, { being fixed points of $\rho(\alpha\beta)$, must be
contained in $\MM=L_\rho$. Thus $\alpha\beta$ generates the
stabilizer of $I_1$.}
Since $\rho$ is a slight perturbation of $\rho_M$, we can assume that
$\sigma_\pm$ are very near to $e$ (in Figure 3). 
%The axis of
%$\rho(\alpha\beta)$ is a geodesic joining $\sigma_{\pm}$ and is depicted
%by the dotted line in figure 3.

}

\begin{figure}[h]
%WinTpicVersion4.30a
{\unitlength 0.1in%
\begin{picture}( 27.9000, 29.8700)( 16.5000,-28.7000)%
% DOT 2 0 3 0 Black White
% 1 3146 1356
% 
\special{pn 4}%
\special{sh 1}%
\special{ar 3146 1356 8 8 0  6.28318530717959E+0000}%
% DOT 0 0 3 0 Black White
% 2 3146 1356 3146 1356
% 
\special{pn 4}%
\special{sh 1}%
\special{ar 3146 1356 16 16 0  6.28318530717959E+0000}%
\special{sh 1}%
\special{ar 3146 1356 16 16 0  6.28318530717959E+0000}%
% CIRCLE 2 0 3 0 Black White
% 4 3146 1356 3146 1356 3140 271 3159 264
% 
\special{pn 8}%
\special{ar 3146 1356 0 0  4.7242932  4.7068591}%
% CIRCLE 2 0 3 0 Black White
% 4 3146 1356 4154 713 4154 713 4154 713
% 
\special{pn 8}%
\special{ar 3146 1356 1196 1196  0.0000000  6.2831853}%
% LINE 2 0 3 0 Black White
% 2 3146 160 3146 2559
% 
\special{pn 8}%
\special{pa 3146 160}%
\special{pa 3146 2559}%
\special{fp}%
% DOT 0 0 3 0 Black White
% 2 2613 1350 2613 1350
% 
\special{pn 4}%
\special{sh 1}%
\special{ar 2613 1350 16 16 0  6.28318530717959E+0000}%
\special{sh 1}%
\special{ar 2613 1350 16 16 0  6.28318530717959E+0000}%
% DOT 2 0 3 0 Black White
% 2 3686 1356 3686 1356
% 
\special{pn 4}%
\special{sh 1}%
\special{ar 3686 1356 8 8 0  6.28318530717959E+0000}%
\special{sh 1}%
\special{ar 3686 1356 8 8 0  6.28318530717959E+0000}%
% SPLINE 2 0 3 0 Black White
% 4 3146 2552 2925 1779 2626 1356 2626 1356
% 
\special{pn 8}%
\special{pa 3146 2552}%
\special{pa 3128 2456}%
\special{pa 3122 2423}%
\special{pa 3110 2359}%
\special{pa 3103 2328}%
\special{pa 3097 2296}%
\special{pa 3083 2232}%
\special{pa 3075 2201}%
\special{pa 3068 2170}%
\special{pa 3060 2139}%
\special{pa 3033 2046}%
\special{pa 3024 2016}%
\special{pa 3013 1986}%
\special{pa 3003 1956}%
\special{pa 2992 1926}%
\special{pa 2968 1868}%
\special{pa 2955 1839}%
\special{pa 2941 1811}%
\special{pa 2927 1782}%
\special{pa 2912 1755}%
\special{pa 2896 1727}%
\special{pa 2880 1700}%
\special{pa 2846 1646}%
\special{pa 2810 1594}%
\special{pa 2753 1516}%
\special{pa 2733 1490}%
\special{pa 2693 1440}%
\special{pa 2673 1414}%
\special{pa 2653 1389}%
\special{pa 2632 1364}%
\special{pa 2626 1356}%
\special{fp}%
% SPLINE 2 0 3 0 Black White
% 4 2620 1356 2912 680 2977 180 2977 173
% 
\special{pn 8}%
\special{pa 2620 1356}%
\special{pa 2635 1327}%
\special{pa 2650 1299}%
\special{pa 2664 1270}%
\special{pa 2679 1241}%
\special{pa 2694 1213}%
\special{pa 2708 1184}%
\special{pa 2723 1155}%
\special{pa 2751 1097}%
\special{pa 2764 1068}%
\special{pa 2778 1039}%
\special{pa 2804 981}%
\special{pa 2817 951}%
\special{pa 2829 922}%
\special{pa 2841 892}%
\special{pa 2874 802}%
\special{pa 2894 742}%
\special{pa 2912 680}%
\special{pa 2928 618}%
\special{pa 2935 587}%
\special{pa 2941 556}%
\special{pa 2947 524}%
\special{pa 2952 492}%
\special{pa 2957 461}%
\special{pa 2965 397}%
\special{pa 2971 333}%
\special{pa 2975 269}%
\special{pa 2977 205}%
\special{pa 2977 173}%
\special{fp}%
% LINE 2 0 3 0 Black White
% 2 2613 1356 1963 1356
% 
\special{pn 8}%
\special{pa 2613 1356}%
\special{pa 1963 1356}%
\special{fp}%
% SPLINE 2 0 3 0 Black White
% 4 1963 1363 2522 1733 2971 2539 2971 2539
% 
\special{pn 8}%
\special{pa 1963 1363}%
\special{pa 1992 1379}%
\special{pa 2021 1394}%
\special{pa 2050 1410}%
\special{pa 2079 1425}%
\special{pa 2107 1441}%
\special{pa 2136 1457}%
\special{pa 2192 1489}%
\special{pa 2248 1523}%
\special{pa 2302 1557}%
\special{pa 2329 1575}%
\special{pa 2355 1593}%
\special{pa 2380 1612}%
\special{pa 2406 1631}%
\special{pa 2478 1691}%
\special{pa 2501 1712}%
\special{pa 2545 1756}%
\special{pa 2566 1779}%
\special{pa 2606 1827}%
\special{pa 2644 1877}%
\special{pa 2662 1902}%
\special{pa 2680 1928}%
\special{pa 2697 1955}%
\special{pa 2713 1982}%
\special{pa 2730 2009}%
\special{pa 2745 2037}%
\special{pa 2761 2065}%
\special{pa 2776 2093}%
\special{pa 2791 2122}%
\special{pa 2819 2180}%
\special{pa 2847 2240}%
\special{pa 2886 2330}%
\special{pa 2899 2361}%
\special{pa 2912 2391}%
\special{pa 2924 2422}%
\special{pa 2937 2453}%
\special{pa 2961 2515}%
\special{pa 2971 2539}%
\special{fp}%
% SPLINE 2 0 3 0 Black White
% 4 3146 173 3406 973 3679 1356 3679 1356
% 
\special{pn 8}%
\special{pa 3146 173}%
\special{pa 3160 237}%
\special{pa 3167 268}%
\special{pa 3175 300}%
\special{pa 3182 332}%
\special{pa 3190 364}%
\special{pa 3197 395}%
\special{pa 3205 427}%
\special{pa 3221 489}%
\special{pa 3230 520}%
\special{pa 3238 552}%
\special{pa 3247 582}%
\special{pa 3257 613}%
\special{pa 3266 644}%
\special{pa 3276 674}%
\special{pa 3287 704}%
\special{pa 3297 734}%
\special{pa 3309 764}%
\special{pa 3320 793}%
\special{pa 3333 823}%
\special{pa 3345 852}%
\special{pa 3358 881}%
\special{pa 3372 909}%
\special{pa 3417 993}%
\special{pa 3433 1020}%
\special{pa 3467 1074}%
\special{pa 3503 1126}%
\special{pa 3560 1204}%
\special{pa 3580 1230}%
\special{pa 3599 1255}%
\special{pa 3619 1280}%
\special{pa 3639 1306}%
\special{pa 3679 1356}%
\special{fp}%
% SPLINE 2 0 3 0 Black White
% 4 3686 1369 3380 1850 3302 2533 3302 2533
% 
\special{pn 8}%
\special{pa 3686 1369}%
\special{pa 3586 1499}%
\special{pa 3567 1525}%
\special{pa 3548 1552}%
\special{pa 3529 1578}%
\special{pa 3493 1632}%
\special{pa 3476 1659}%
\special{pa 3444 1713}%
\special{pa 3429 1741}%
\special{pa 3415 1769}%
\special{pa 3402 1798}%
\special{pa 3389 1826}%
\special{pa 3378 1855}%
\special{pa 3358 1915}%
\special{pa 3342 1975}%
\special{pa 3324 2068}%
\special{pa 3320 2099}%
\special{pa 3316 2131}%
\special{pa 3310 2195}%
\special{pa 3308 2227}%
\special{pa 3306 2260}%
\special{pa 3305 2293}%
\special{pa 3304 2325}%
\special{pa 3302 2391}%
\special{pa 3302 2533}%
\special{fp}%
% LINE 2 0 3 0 Black White
% 2 3679 1356 4342 1363
% 
\special{pn 8}%
\special{pa 3679 1356}%
\special{pa 4342 1363}%
\special{fp}%
% SPLINE 2 0 3 0 Black White
% 4 3309 2533 3653 1857 4323 1551 4323 1551
% 
\special{pn 8}%
\special{pa 3309 2533}%
\special{pa 3320 2501}%
\special{pa 3331 2468}%
\special{pa 3375 2340}%
\special{pa 3387 2309}%
\special{pa 3399 2277}%
\special{pa 3411 2247}%
\special{pa 3424 2216}%
\special{pa 3450 2156}%
\special{pa 3464 2127}%
\special{pa 3478 2099}%
\special{pa 3493 2070}%
\special{pa 3525 2016}%
\special{pa 3542 1990}%
\special{pa 3578 1940}%
\special{pa 3597 1917}%
\special{pa 3617 1894}%
\special{pa 3638 1872}%
\special{pa 3660 1851}%
\special{pa 3682 1831}%
\special{pa 3730 1793}%
\special{pa 3755 1776}%
\special{pa 3781 1760}%
\special{pa 3808 1744}%
\special{pa 3835 1729}%
\special{pa 3863 1714}%
\special{pa 3892 1700}%
\special{pa 3921 1687}%
\special{pa 3951 1674}%
\special{pa 4011 1650}%
\special{pa 4042 1639}%
\special{pa 4074 1628}%
\special{pa 4105 1617}%
\special{pa 4169 1597}%
\special{pa 4202 1587}%
\special{pa 4234 1577}%
\special{pa 4267 1567}%
\special{pa 4299 1558}%
\special{pa 4323 1551}%
\special{fp}%
% SPLINE 2 2 3 0 Black White
% 4 2561 2403 3140 2169 3751 2383 3751 2383
% 
\special{pn 8}%
\special{pn 8}%
\special{pa 2561 2403}%
\special{pa 2568 2399}%
\special{fp}%
\special{pa 2600 2380}%
\special{pa 2607 2376}%
\special{fp}%
\special{pa 2640 2357}%
\special{pa 2647 2353}%
\special{fp}%
\special{pa 2679 2334}%
\special{pa 2686 2330}%
\special{fp}%
\special{pa 2719 2311}%
\special{pa 2726 2308}%
\special{fp}%
\special{pa 2759 2290}%
\special{pa 2766 2286}%
\special{fp}%
\special{pa 2800 2269}%
\special{pa 2807 2265}%
\special{fp}%
\special{pa 2841 2249}%
\special{pa 2848 2246}%
\special{fp}%
\special{pa 2883 2231}%
\special{pa 2890 2228}%
\special{fp}%
\special{pa 2925 2215}%
\special{pa 2933 2212}%
\special{fp}%
\special{pa 2968 2199}%
\special{pa 2976 2197}%
\special{fp}%
\special{pa 3012 2187}%
\special{pa 3020 2185}%
\special{fp}%
\special{pa 3057 2178}%
\special{pa 3065 2176}%
\special{fp}%
\special{pa 3102 2172}%
\special{pa 3110 2171}%
\special{fp}%
\special{pa 3148 2169}%
\special{pa 3156 2169}%
\special{fp}%
\special{pa 3193 2170}%
\special{pa 3201 2171}%
\special{fp}%
\special{pa 3239 2175}%
\special{pa 3247 2176}%
\special{fp}%
\special{pa 3284 2183}%
\special{pa 3292 2184}%
\special{fp}%
\special{pa 3328 2193}%
\special{pa 3336 2195}%
\special{fp}%
\special{pa 3372 2206}%
\special{pa 3380 2209}%
\special{fp}%
\special{pa 3415 2221}%
\special{pa 3423 2224}%
\special{fp}%
\special{pa 3458 2238}%
\special{pa 3465 2241}%
\special{fp}%
\special{pa 3499 2256}%
\special{pa 3507 2259}%
\special{fp}%
\special{pa 3541 2275}%
\special{pa 3548 2278}%
\special{fp}%
\special{pa 3582 2295}%
\special{pa 3589 2298}%
\special{fp}%
\special{pa 3623 2315}%
\special{pa 3630 2319}%
\special{fp}%
\special{pa 3664 2336}%
\special{pa 3671 2340}%
\special{fp}%
\special{pa 3703 2358}%
\special{pa 3711 2362}%
\special{fp}%
\special{pa 3744 2379}%
\special{pa 3751 2383}%
\special{fp}%
% CIRCLE 2 0 3 0 Black White
% 4 3146 1356 4394 953 4394 953 3978 329
% 
\special{pn 8}%
\special{ar 3146 1356 1311 1311  5.3932742  5.9708389}%
% VECTOR 2 0 3 0 Black White
% 2 4394 953 4420 1070
% 
\special{pn 8}%
\special{pa 4394 953}%
\special{pa 4420 1070}%
\special{fp}%
\special{sh 1}%
\special{pa 4420 1070}%
\special{pa 4425 1001}%
\special{pa 4408 1018}%
\special{pa 4386 1009}%
\special{pa 4420 1070}%
\special{fp}%
% STR 2 0 3 0 Black White
% 4 2860 1250 2860 1350 2 0 0 0
% $P$
\put(28.6000,-13.5000){\makebox(0,0)[lb]{$P$}}%
% STR 2 0 3 0 Black White
% 4 2400 2400 2400 2500 2 0 0 0
% $\sigma_-$
\put(24.0000,-25.0000){\makebox(0,0)[lb]{$\sigma_-$}}%
% STR 2 0 3 0 Black White
% 4 2790 2600 2790 2700 2 0 0 0
% $\beta^2\alpha$
\put(27.9000,-27.0000){\makebox(0,0)[lb]{$\beta^2\alpha$}}%
% STR 2 0 3 0 Black White
% 4 3710 2420 3710 2520 2 0 0 0
% $\sigma_+$
\put(37.1000,-25.2000){\makebox(0,0)[lb]{$\sigma_+$}}%
% STR 2 0 3 0 Black White
% 4 4420 1520 4420 1620 2 0 0 0
% $\alpha\beta\alpha$
\put(44.2000,-16.2000){\makebox(0,0)[lb]{$\alpha\beta\alpha$}}%
% STR 2 0 3 0 Black White
% 4 4440 1310 4440 1410 2 0 0 0
% $\alpha\beta^2$
\put(44.4000,-14.1000){\makebox(0,0)[lb]{$\alpha\beta^2$}}%
% STR 2 0 3 0 Black White
% 4 3230 1320 3230 1420 2 0 0 0
% $i$
\put(32.3000,-14.2000){\makebox(0,0)[lb]{$i$}}%
% STR 2 0 3 0 Black White
% 4 2510 1190 2510 1290 2 0 0 0
% $\omega'$
\put(25.1000,-12.9000){\makebox(0,0)[lb]{$\omega'$}}%
% STR 2 0 3 0 Black White
% 4 3150 0 3150 100 2 0 0 0
% $\alpha$
\put(31.5000,-1.0000){\makebox(0,0)[lb]{$\alpha$}}%
% STR 2 0 3 0 Black White
% 4 2830 10 2830 110 2 0 0 0
% $\beta$
\put(28.3000,-1.1000){\makebox(0,0)[lb]{$\beta$}}%
% STR 2 0 3 0 Black White
% 4 1710 1340 1710 1440 2 0 0 0
% $\beta^2$
\put(17.1000,-14.4000){\makebox(0,0)[lb]{$\beta^2$}}%
% STR 2 0 3 0 Black White
% 4 3400 1290 3400 1390 2 0 0 0
% $\alpha P$
\put(34.0000,-13.9000){\makebox(0,0)[lb]{$\alpha P$}}%
% STR 2 0 3 0 Black White
% 4 3590 1620 3590 1720 2 0 0 0
% $\alpha\beta P$
\put(35.9000,-17.2000){\makebox(0,0)[lb]{$\alpha\beta P$}}%
% STR 2 0 3 0 Black White
% 4 2600 1640 2600 1740 2 0 0 0
% $\beta^2 P$
\put(26.0000,-17.4000){\makebox(0,0)[lb]{$\beta^2 P$}}%
% STR 2 0 3 0 Black White
% 4 3110 2590 3110 2690 2 0 0 0
% $e$
\put(31.1000,-26.9000){\makebox(0,0)[lb]{$e$}}%
% STR 2 0 3 0 Black White
% 4 3280 2560 3280 2660 2 0 0 0
% $\alpha\beta$
\put(32.8000,-26.6000){\makebox(0,0)[lb]{$\alpha\beta$}}%
% CIRCLE 2 0 3 0 Black White
% 4 3140 1360 2670 2760 2690 2780 3630 2800
% 
\special{pn 8}%
\special{ar 3140 1360 1477 1477  1.2428088  1.8776860}%
% VECTOR 2 0 3 0 Black White
% 2 3600 2760 3690 2740
% 
\special{pn 8}%
\special{pa 3600 2760}%
\special{pa 3690 2740}%
\special{fp}%
\special{sh 1}%
\special{pa 3690 2740}%
\special{pa 3621 2735}%
\special{pa 3638 2752}%
\special{pa 3629 2774}%
\special{pa 3690 2740}%
\special{fp}%
% SPLINE 2 0 3 0 Black White
% 4 2970 180 2590 880 1970 1160 1970 1160
% 
\special{pn 8}%
\special{pa 2970 180}%
\special{pa 2958 212}%
\special{pa 2947 244}%
\special{pa 2935 276}%
\special{pa 2924 308}%
\special{pa 2900 372}%
\special{pa 2887 403}%
\special{pa 2875 434}%
\special{pa 2862 465}%
\special{pa 2849 495}%
\special{pa 2821 555}%
\special{pa 2807 584}%
\special{pa 2777 640}%
\special{pa 2761 668}%
\special{pa 2744 695}%
\special{pa 2727 721}%
\special{pa 2710 746}%
\special{pa 2691 771}%
\special{pa 2672 795}%
\special{pa 2652 818}%
\special{pa 2610 862}%
\special{pa 2564 902}%
\special{pa 2540 920}%
\special{pa 2515 938}%
\special{pa 2489 955}%
\special{pa 2463 971}%
\special{pa 2436 986}%
\special{pa 2408 1001}%
\special{pa 2379 1015}%
\special{pa 2351 1028}%
\special{pa 2321 1041}%
\special{pa 2261 1065}%
\special{pa 2230 1077}%
\special{pa 2199 1088}%
\special{pa 2167 1098}%
\special{pa 2136 1109}%
\special{pa 2040 1139}%
\special{pa 2007 1149}%
\special{pa 1975 1159}%
\special{pa 1970 1160}%
\special{fp}%
% SPLINE 2 0 3 0 Black White
% 4 4350 1360 3570 800 3370 180 3370 180
% 
\special{pn 8}%
\special{pa 4350 1360}%
\special{pa 4288 1330}%
\special{pa 4258 1314}%
\special{pa 4227 1299}%
\special{pa 4197 1284}%
\special{pa 4166 1268}%
\special{pa 4136 1253}%
\special{pa 4106 1237}%
\special{pa 4077 1221}%
\special{pa 4047 1205}%
\special{pa 4018 1188}%
\special{pa 3990 1171}%
\special{pa 3961 1154}%
\special{pa 3933 1137}%
\special{pa 3906 1120}%
\special{pa 3879 1102}%
\special{pa 3852 1083}%
\special{pa 3827 1065}%
\special{pa 3801 1045}%
\special{pa 3777 1026}%
\special{pa 3753 1006}%
\special{pa 3729 985}%
\special{pa 3706 964}%
\special{pa 3684 942}%
\special{pa 3663 920}%
\special{pa 3643 898}%
\special{pa 3623 874}%
\special{pa 3587 826}%
\special{pa 3570 800}%
\special{pa 3554 774}%
\special{pa 3539 747}%
\special{pa 3525 720}%
\special{pa 3499 664}%
\special{pa 3488 635}%
\special{pa 3476 605}%
\special{pa 3456 545}%
\special{pa 3438 483}%
\special{pa 3429 451}%
\special{pa 3421 420}%
\special{pa 3414 388}%
\special{pa 3406 355}%
\special{pa 3399 323}%
\special{pa 3392 290}%
\special{pa 3386 258}%
\special{pa 3372 192}%
\special{pa 3370 180}%
\special{fp}%
% STR 2 0 3 0 Black White
% 4 1650 1120 1650 1220 2 0 0 0
% $\beta\alpha$
\put(16.5000,-12.2000){\makebox(0,0)[lb]{$\beta\alpha$}}%
% STR 2 0 3 0 Black White
% 4 3430 -10 3430 90 2 0 0 0
% $\alpha\beta^2\alpha$
\put(34.3000,-0.9000){\makebox(0,0)[lb]{$\alpha\beta^2\alpha$}}%
% STR 2 0 3 0 Black White
% 4 4350 540 4350 640 2 0 0 0
% $\rho(\alpha\beta)$
\put(43.5000,-6.4000){\makebox(0,0)[lb]{$\rho(\alpha\beta)$}}%
% STR 2 0 3 0 Black White
% 4 2920 2900 2920 3000 2 0 0 0
% $\rho(\alpha\beta)$
\put(29.2000,-30.0000){\makebox(0,0)[lb]{$\rho(\alpha\beta)$}}%
\end{picture}}%
\caption{$e$ indicates the base point $y_0=0$ and $g$ indicates
the point $\rho(g)y_0$. $P$ is a fundamental domain of
 $\rho(\overline G)$. }
\end{figure}
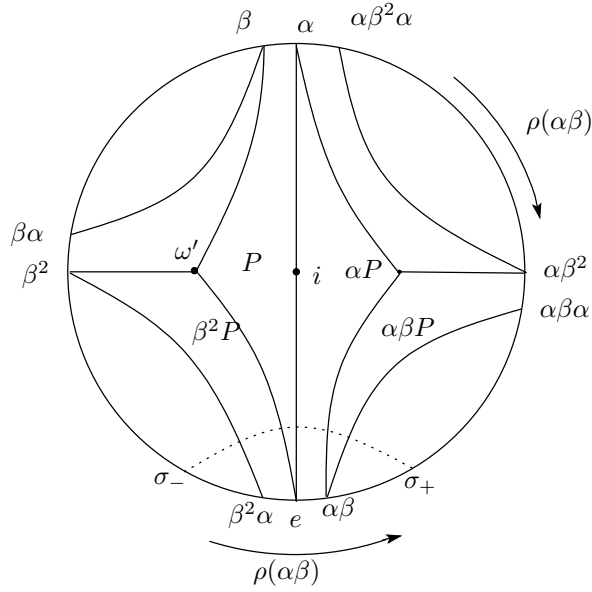

Since the translates of $P\vert_{S^1}$ tesselate $S^1\setminus\MM$,
there are no gaps of $\MM$ other than the translates of $I_1$.
{  In particular, any gap of $\Cl(\rho(\oG)y_0)$ is a gap of
$\rho(\oG)y_0$. 
This, together with the freeness of the action $\rho$ at $y_0$,
 shows that $\rho$ acts tightly at
$y_0$.} Thus $\rho$ is topologically conjugate to
a dynamical realization of a circular order
$c\in CO(\oG)$. The linear part of $c$ is the subgroup 
$\langle\alpha\beta\rangle$ by Lemma \ref{linear-part}.

\medskip
Let us show first of all that $c$ is isolated.
In \cite{MR}, Proposition 3.3, the authors showed that the dynamical
realization is continuous. More precisely, they showed the following.

\begin{proposition}\label{continuous}
Given any neighbourhood $U$ of $\rho$ in ${\rm Hom}(\oG,{\rm  Homeo}_+(S^1))$,
there is a neighbourhood $V$ of $c$ in $CO(\oG)$ such that any element
in $V$ has a conjugate of its dynamical realization contained in $U$.
\end{proposition}

{ 
 In Figure 3, six translates of $P$ are depicted. Let 
$Q$ be their union:
$$
Q=P\cup \beta P\cup\beta^2P\cup\alpha P\cup\alpha\beta P
\cup\alpha\beta^2P.$$
The convex set $Q\cap\mathbb H$ has four sides. 
Let $\gamma_1=\beta^2\alpha\beta\alpha$ and
$\gamma_2=\alpha\beta\alpha\beta^2$.
Then $\rho(\gamma_i)$ maps a side of $Q\cap\mathbb H$ onto its opposite side,
as indicated in Figure 4. By the Klein criterion (also known as the 
ping-pong lemma), $\gamma_1$ and $\gamma_2$ are free generators of a free subgroup  $H$
of $\oG$. Moreover, for any nontrivial element $\gamma$ of $H$, 
$\rho(\gamma)$ is a hyperbolic transformation.
The subgroup $H$ coincides with the commutator subgroup $[\oG,
\oG]$.
{ To show this, notice that $\gamma_i$ are commutators and thus
$H\subset[\oG,\oG]$. On the other hand, $H$ is a index 6 subgroup of
$\oG$, since the fundamental domain $Q$ of $\rho(H)$ consists of 6 iterates of $P$,
and $[\oG,\oG]$ is also of index 6, since the abelianization of $\oG$ is 
$\Z/6\Z$.
}

Below we indicate the point $\rho(g)y_0$ simply by $g$, as we already
did in the figures. Let us define four intervals of $S^1$:
$$
K_1^{-}=[\alpha\beta^2,\alpha\beta^2\alpha],\ \
K_1^+=[\beta^2,\beta^2\alpha], \ \
K_2^-=[\beta,\beta\alpha],\ \ K_2^+=[\alpha\beta,\alpha\beta\alpha].$$
Then $\rho(\gamma_1)(S^1\setminus K_1^-)=\Int(K_1^+)$ and
$\rho(\gamma_2)(S^1\setminus K_2^-)=\Int(K_2^+)$. 
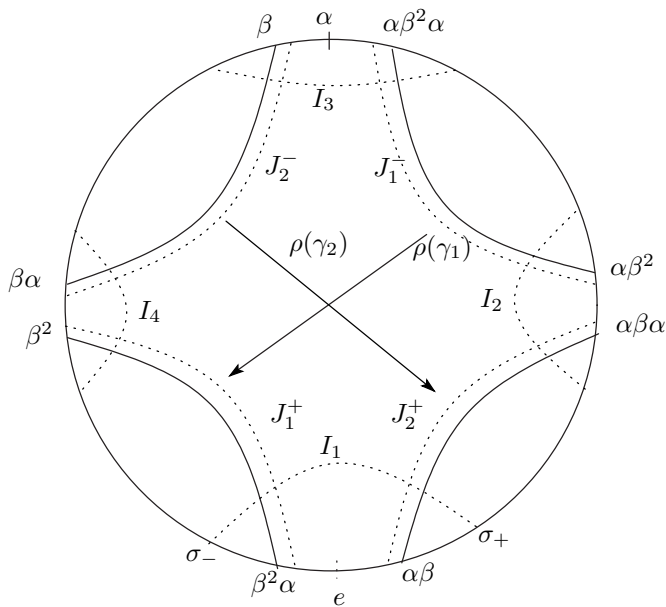
\begin{figure}[h]
%WinTpicVersion4.30a
{\unitlength 0.1in%
\begin{picture}( 31.8000, 30.4000)( 50.2000,-35.7000)%
% CIRCLE 2 0 3 0 Black White
% 4 6710 2130 7780 1240 7780 1240 7780 1240
% 
\special{pn 8}%
\special{ar 6710 2130 1392 1392  0.0000000  6.2831853}%
% SPLINE 2 0 3 0 Black White
% 4 6420 770 6040 1660 5330 2020 5330 2020
% 
\special{pn 8}%
\special{pa 6420 770}%
\special{pa 6412 804}%
\special{pa 6404 837}%
\special{pa 6395 870}%
\special{pa 6387 904}%
\special{pa 6378 937}%
\special{pa 6370 970}%
\special{pa 6343 1069}%
\special{pa 6333 1101}%
\special{pa 6324 1133}%
\special{pa 6314 1165}%
\special{pa 6281 1258}%
\special{pa 6257 1318}%
\special{pa 6245 1347}%
\special{pa 6232 1376}%
\special{pa 6218 1405}%
\special{pa 6204 1433}%
\special{pa 6189 1460}%
\special{pa 6173 1487}%
\special{pa 6157 1513}%
\special{pa 6141 1538}%
\special{pa 6123 1563}%
\special{pa 6105 1587}%
\special{pa 6086 1610}%
\special{pa 6066 1633}%
\special{pa 6045 1655}%
\special{pa 6024 1676}%
\special{pa 5978 1716}%
\special{pa 5954 1734}%
\special{pa 5904 1770}%
\special{pa 5878 1787}%
\special{pa 5851 1803}%
\special{pa 5824 1818}%
\special{pa 5796 1833}%
\special{pa 5767 1848}%
\special{pa 5739 1862}%
\special{pa 5709 1876}%
\special{pa 5619 1915}%
\special{pa 5557 1939}%
\special{pa 5525 1951}%
\special{pa 5494 1962}%
\special{pa 5462 1974}%
\special{pa 5398 1996}%
\special{pa 5366 2008}%
\special{pa 5334 2019}%
\special{pa 5330 2020}%
\special{fp}%
% SPLINE 2 0 3 0 Black White
% 4 7030 790 7330 1630 8090 1960 8090 1960
% 
\special{pn 8}%
\special{pa 7030 790}%
\special{pa 7036 825}%
\special{pa 7042 859}%
\special{pa 7048 894}%
\special{pa 7060 962}%
\special{pa 7067 997}%
\special{pa 7074 1031}%
\special{pa 7081 1064}%
\special{pa 7088 1098}%
\special{pa 7096 1131}%
\special{pa 7104 1163}%
\special{pa 7112 1196}%
\special{pa 7121 1228}%
\special{pa 7130 1259}%
\special{pa 7140 1290}%
\special{pa 7151 1321}%
\special{pa 7162 1350}%
\special{pa 7173 1380}%
\special{pa 7185 1408}%
\special{pa 7198 1436}%
\special{pa 7212 1464}%
\special{pa 7242 1516}%
\special{pa 7258 1541}%
\special{pa 7294 1589}%
\special{pa 7313 1611}%
\special{pa 7333 1633}%
\special{pa 7354 1653}%
\special{pa 7376 1673}%
\special{pa 7399 1691}%
\special{pa 7423 1709}%
\special{pa 7447 1726}%
\special{pa 7473 1743}%
\special{pa 7499 1758}%
\special{pa 7526 1773}%
\special{pa 7554 1787}%
\special{pa 7583 1801}%
\special{pa 7612 1814}%
\special{pa 7672 1838}%
\special{pa 7734 1860}%
\special{pa 7766 1871}%
\special{pa 7830 1891}%
\special{pa 7896 1909}%
\special{pa 7929 1919}%
\special{pa 7963 1927}%
\special{pa 7997 1936}%
\special{pa 8030 1945}%
\special{pa 8064 1953}%
\special{pa 8090 1960}%
\special{fp}%
% SPLINE 2 0 3 0 Black White
% 4 5330 2300 6070 2640 6430 3510 6430 3510
% 
\special{pn 8}%
\special{pa 5330 2300}%
\special{pa 5363 2310}%
\special{pa 5396 2319}%
\special{pa 5462 2339}%
\special{pa 5494 2349}%
\special{pa 5527 2359}%
\special{pa 5559 2369}%
\special{pa 5591 2380}%
\special{pa 5622 2391}%
\special{pa 5654 2402}%
\special{pa 5685 2413}%
\special{pa 5715 2425}%
\special{pa 5746 2438}%
\special{pa 5775 2450}%
\special{pa 5804 2464}%
\special{pa 5833 2477}%
\special{pa 5889 2507}%
\special{pa 5915 2522}%
\special{pa 5941 2538}%
\special{pa 5967 2555}%
\special{pa 6015 2591}%
\special{pa 6038 2611}%
\special{pa 6060 2631}%
\special{pa 6082 2652}%
\special{pa 6102 2673}%
\special{pa 6140 2719}%
\special{pa 6158 2744}%
\special{pa 6175 2769}%
\special{pa 6191 2794}%
\special{pa 6207 2821}%
\special{pa 6221 2848}%
\special{pa 6236 2875}%
\special{pa 6275 2962}%
\special{pa 6287 2992}%
\special{pa 6298 3022}%
\special{pa 6320 3084}%
\special{pa 6330 3116}%
\special{pa 6339 3148}%
\special{pa 6349 3180}%
\special{pa 6358 3213}%
\special{pa 6367 3245}%
\special{pa 6375 3279}%
\special{pa 6384 3312}%
\special{pa 6392 3345}%
\special{pa 6416 3447}%
\special{pa 6423 3480}%
\special{pa 6430 3510}%
\special{fp}%
% SPLINE 2 0 3 0 Black White
% 4 8110 2280 7420 2700 7080 3480 7080 3480
% 
\special{pn 8}%
\special{pa 8110 2280}%
\special{pa 8079 2294}%
\special{pa 8048 2307}%
\special{pa 8018 2321}%
\special{pa 7987 2334}%
\special{pa 7956 2348}%
\special{pa 7896 2376}%
\special{pa 7866 2391}%
\special{pa 7836 2405}%
\special{pa 7778 2435}%
\special{pa 7749 2451}%
\special{pa 7693 2483}%
\special{pa 7612 2534}%
\special{pa 7587 2552}%
\special{pa 7561 2571}%
\special{pa 7513 2609}%
\special{pa 7467 2651}%
\special{pa 7446 2672}%
\special{pa 7425 2694}%
\special{pa 7405 2717}%
\special{pa 7386 2741}%
\special{pa 7367 2766}%
\special{pa 7333 2816}%
\special{pa 7316 2843}%
\special{pa 7286 2897}%
\special{pa 7271 2925}%
\special{pa 7258 2953}%
\special{pa 7244 2982}%
\special{pa 7232 3012}%
\special{pa 7219 3041}%
\special{pa 7207 3071}%
\special{pa 7174 3164}%
\special{pa 7164 3195}%
\special{pa 7153 3227}%
\special{pa 7144 3259}%
\special{pa 7124 3323}%
\special{pa 7115 3355}%
\special{pa 7106 3388}%
\special{pa 7097 3420}%
\special{pa 7088 3453}%
\special{pa 7080 3480}%
\special{fp}%
% LINE 2 0 3 0 Black White
% 2 6700 700 6700 790
% 
\special{pn 8}%
\special{pa 6700 700}%
\special{pa 6700 790}%
\special{fp}%
% LINE 2 2 3 0 Black White
% 2 6740 3470 6740 3560
% 
\special{pn 8}%
\special{pa 6740 3470}%
\special{pa 6740 3560}%
\special{dt 0.045}%
% VECTOR 2 0 3 0 Black White
% 2 7210 1760 6180 2500
% 
\special{pn 8}%
\special{pa 7210 1760}%
\special{pa 6180 2500}%
\special{fp}%
\special{sh 1}%
\special{pa 6180 2500}%
\special{pa 6246 2477}%
\special{pa 6223 2469}%
\special{pa 6222 2445}%
\special{pa 6180 2500}%
\special{fp}%
% VECTOR 2 0 3 0 Black White
% 2 6160 1690 7250 2580
% 
\special{pn 8}%
\special{pa 6160 1690}%
\special{pa 7250 2580}%
\special{fp}%
\special{sh 1}%
\special{pa 7250 2580}%
\special{pa 7211 2522}%
\special{pa 7209 2546}%
\special{pa 7186 2553}%
\special{pa 7250 2580}%
\special{fp}%
% STR 2 0 3 0 Black White
% 4 6860 1800 6860 1900 2 0 0 0
% $\rho(\gamma_1)$
\put(68.6000,-19.0000){\makebox(0,0)[lb]{$\ \ \ \ \ \ \rho(\gamma_1)$}}%
% STR 2 0 3 0 Black White
% 4 6490 1780 6490 1880 2 0 0 0
% $\rho(\gamma_2)$
\put(64.9000,-18.8000){\makebox(0,0)[lb]{$\rho(\gamma_2)$}}%
% STR 2 0 3 0 Black White
% 4 6320 630 6320 730 2 0 0 0
% $\beta$
\put(63.2000,-7.3000){\makebox(0,0)[lb]{$\beta$}}%
% STR 2 0 3 0 Black White
% 4 6630 560 6630 660 2 0 0 0
% $\alpha$
\put(66.3000,-6.6000){\makebox(0,0)[lb]{$\alpha$}}%
% STR 2 0 3 0 Black White
% 4 8170 1900 8170 2000 2 0 0 0
% $\alpha\beta^2$
\put(81.7000,-20.0000){\makebox(0,0)[lb]{$\alpha\beta^2$}}%
% STR 2 0 3 0 Black White
% 4 8200 2200 8200 2300 2 0 0 0
% $\alpha\beta\alpha$
\put(82.0000,-23.0000){\makebox(0,0)[lb]{$\alpha\beta\alpha$}}%
% STR 2 0 3 0 Black White
% 4 7080 3500 7080 3600 2 0 0 0
% $\alpha\beta$
\put(70.8000,-36.0000){\makebox(0,0)[lb]{$\alpha\beta$}}%
% STR 2 0 3 0 Black White
% 4 6720 3600 6720 3700 2 0 0 0
% $e$
\put(67.2000,-37.0000){\makebox(0,0)[lb]{$e$}}%
% STR 2 0 3 0 Black White
% 4 5110 2260 5110 2360 2 0 0 0
% $\beta^2$
\put(51.1000,-23.6000){\makebox(0,0)[lb]{$\beta^2$}}%
% STR 2 0 3 0 Black White
% 4 5020 1960 5020 2060 2 0 0 0
% $\beta\alpha$
\put(50.2000,-20.6000){\makebox(0,0)[lb]{$\beta\alpha$}}%
% SPLINE 2 2 3 0 Black White
% 4 6070 3380 6720 2960 7470 3290 7470 3290
% 
\special{pn 8}%
\special{pn 8}%
\special{pa 6070 3380}%
\special{pa 6076 3374}%
\special{fp}%
\special{pa 6103 3349}%
\special{pa 6109 3344}%
\special{fp}%
\special{pa 6136 3319}%
\special{pa 6141 3313}%
\special{fp}%
\special{pa 6168 3288}%
\special{pa 6174 3283}%
\special{fp}%
\special{pa 6202 3258}%
\special{pa 6208 3253}%
\special{fp}%
\special{pa 6235 3228}%
\special{pa 6241 3223}%
\special{fp}%
\special{pa 6269 3199}%
\special{pa 6275 3194}%
\special{fp}%
\special{pa 6303 3170}%
\special{pa 6310 3165}%
\special{fp}%
\special{pa 6338 3142}%
\special{pa 6345 3137}%
\special{fp}%
\special{pa 6374 3115}%
\special{pa 6380 3110}%
\special{fp}%
\special{pa 6410 3089}%
\special{pa 6417 3084}%
\special{fp}%
\special{pa 6448 3064}%
\special{pa 6455 3060}%
\special{fp}%
\special{pa 6487 3041}%
\special{pa 6493 3037}%
\special{fp}%
\special{pa 6526 3020}%
\special{pa 6533 3016}%
\special{fp}%
\special{pa 6567 3001}%
\special{pa 6574 2998}%
\special{fp}%
\special{pa 6608 2985}%
\special{pa 6616 2982}%
\special{fp}%
\special{pa 6651 2972}%
\special{pa 6659 2970}%
\special{fp}%
\special{pa 6695 2963}%
\special{pa 6703 2962}%
\special{fp}%
\special{pa 6740 2958}%
\special{pa 6748 2958}%
\special{fp}%
\special{pa 6785 2958}%
\special{pa 6793 2959}%
\special{fp}%
\special{pa 6829 2962}%
\special{pa 6837 2964}%
\special{fp}%
\special{pa 6874 2971}%
\special{pa 6881 2972}%
\special{fp}%
\special{pa 6917 2981}%
\special{pa 6925 2983}%
\special{fp}%
\special{pa 6960 2995}%
\special{pa 6967 2997}%
\special{fp}%
\special{pa 7002 3010}%
\special{pa 7009 3014}%
\special{fp}%
\special{pa 7043 3028}%
\special{pa 7050 3032}%
\special{fp}%
\special{pa 7084 3047}%
\special{pa 7091 3051}%
\special{fp}%
\special{pa 7123 3068}%
\special{pa 7130 3072}%
\special{fp}%
\special{pa 7163 3090}%
\special{pa 7169 3094}%
\special{fp}%
\special{pa 7202 3112}%
\special{pa 7208 3116}%
\special{fp}%
\special{pa 7239 3136}%
\special{pa 7246 3141}%
\special{fp}%
\special{pa 7278 3160}%
\special{pa 7284 3164}%
\special{fp}%
\special{pa 7315 3184}%
\special{pa 7322 3189}%
\special{fp}%
\special{pa 7353 3209}%
\special{pa 7359 3214}%
\special{fp}%
\special{pa 7389 3235}%
\special{pa 7396 3239}%
\special{fp}%
\special{pa 7426 3260}%
\special{pa 7433 3265}%
\special{fp}%
\special{pa 7463 3286}%
\special{pa 7470 3290}%
\special{fp}%
% STR 2 0 3 0 Black White
% 4 5950 3400 5950 3500 2 0 0 0
% $\sigma_-$
\put(59.5000,-35.0000){\makebox(0,0)[lb]{$\sigma_-$}}%
% STR 2 0 3 0 Black White
% 4 7480 3300 7480 3400 2 0 0 0
% $\sigma_+$
\put(74.8000,-34.0000){\makebox(0,0)[lb]{$\sigma_+$}}%
% STR 2 0 3 0 Black White
% 4 6290 3560 6290 3660 2 0 0 0
% $\beta^2\alpha$
\put(62.9000,-36.6000){\makebox(0,0)[lb]{$\beta^2\alpha$}}%
% STR 2 0 3 0 Black White
% 4 6980 620 6980 720 2 0 0 0
% $\alpha\beta^2\alpha$
\put(69.8000,-7.2000){\makebox(0,0)[lb]{$\alpha\beta^2\alpha$}}%
% SPLINE 2 2 3 0 Black White
% 4 6930 770 7240 1650 8090 2020 8090 2020
% 
\special{pn 8}%
\special{pn 8}%
\special{pa 6930 770}%
\special{pa 6931 778}%
\special{fp}%
\special{pa 6938 815}%
\special{pa 6939 823}%
\special{fp}%
\special{pa 6946 860}%
\special{pa 6947 868}%
\special{fp}%
\special{pa 6954 905}%
\special{pa 6955 913}%
\special{fp}%
\special{pa 6962 950}%
\special{pa 6963 958}%
\special{fp}%
\special{pa 6970 995}%
\special{pa 6972 1003}%
\special{fp}%
\special{pa 6979 1040}%
\special{pa 6981 1048}%
\special{fp}%
\special{pa 6989 1085}%
\special{pa 6991 1093}%
\special{fp}%
\special{pa 7000 1130}%
\special{pa 7002 1138}%
\special{fp}%
\special{pa 7011 1174}%
\special{pa 7013 1182}%
\special{fp}%
\special{pa 7022 1218}%
\special{pa 7024 1226}%
\special{fp}%
\special{pa 7035 1262}%
\special{pa 7038 1270}%
\special{fp}%
\special{pa 7049 1306}%
\special{pa 7052 1314}%
\special{fp}%
\special{pa 7064 1349}%
\special{pa 7067 1357}%
\special{fp}%
\special{pa 7081 1392}%
\special{pa 7084 1399}%
\special{fp}%
\special{pa 7099 1434}%
\special{pa 7103 1441}%
\special{fp}%
\special{pa 7119 1475}%
\special{pa 7123 1482}%
\special{fp}%
\special{pa 7142 1515}%
\special{pa 7146 1522}%
\special{fp}%
\special{pa 7166 1554}%
\special{pa 7170 1561}%
\special{fp}%
\special{pa 7191 1592}%
\special{pa 7196 1598}%
\special{fp}%
\special{pa 7220 1628}%
\special{pa 7225 1634}%
\special{fp}%
\special{pa 7251 1661}%
\special{pa 7256 1667}%
\special{fp}%
\special{pa 7284 1693}%
\special{pa 7290 1698}%
\special{fp}%
\special{pa 7319 1722}%
\special{pa 7326 1727}%
\special{fp}%
\special{pa 7356 1750}%
\special{pa 7363 1754}%
\special{fp}%
\special{pa 7394 1775}%
\special{pa 7401 1779}%
\special{fp}%
\special{pa 7433 1798}%
\special{pa 7441 1802}%
\special{fp}%
\special{pa 7474 1820}%
\special{pa 7481 1824}%
\special{fp}%
\special{pa 7515 1840}%
\special{pa 7522 1843}%
\special{fp}%
\special{pa 7557 1859}%
\special{pa 7564 1862}%
\special{fp}%
\special{pa 7599 1876}%
\special{pa 7607 1879}%
\special{fp}%
\special{pa 7642 1893}%
\special{pa 7649 1895}%
\special{fp}%
\special{pa 7685 1908}%
\special{pa 7693 1911}%
\special{fp}%
\special{pa 7729 1922}%
\special{pa 7736 1924}%
\special{fp}%
\special{pa 7772 1935}%
\special{pa 7780 1938}%
\special{fp}%
\special{pa 7816 1949}%
\special{pa 7824 1951}%
\special{fp}%
\special{pa 7860 1961}%
\special{pa 7868 1963}%
\special{fp}%
\special{pa 7905 1973}%
\special{pa 7912 1975}%
\special{fp}%
\special{pa 7949 1985}%
\special{pa 7957 1987}%
\special{fp}%
\special{pa 7993 1995}%
\special{pa 8001 1998}%
\special{fp}%
\special{pa 8038 2007}%
\special{pa 8046 2009}%
\special{fp}%
\special{pa 8082 2018}%
\special{pa 8090 2020}%
\special{fp}%
% SPLINE 2 2 3 0 Black White
% 4 6500 760 6140 1660 5330 2080 5330 2080
% 
\special{pn 8}%
\special{pn 8}%
\special{pa 6500 760}%
\special{pa 6498 768}%
\special{fp}%
\special{pa 6490 804}%
\special{pa 6488 812}%
\special{fp}%
\special{pa 6480 848}%
\special{pa 6478 856}%
\special{fp}%
\special{pa 6470 892}%
\special{pa 6468 900}%
\special{fp}%
\special{pa 6459 936}%
\special{pa 6457 944}%
\special{fp}%
\special{pa 6448 980}%
\special{pa 6447 987}%
\special{fp}%
\special{pa 6438 1023}%
\special{pa 6436 1031}%
\special{fp}%
\special{pa 6426 1067}%
\special{pa 6424 1075}%
\special{fp}%
\special{pa 6414 1110}%
\special{pa 6412 1118}%
\special{fp}%
\special{pa 6401 1154}%
\special{pa 6398 1161}%
\special{fp}%
\special{pa 6388 1197}%
\special{pa 6385 1204}%
\special{fp}%
\special{pa 6373 1239}%
\special{pa 6371 1247}%
\special{fp}%
\special{pa 6358 1282}%
\special{pa 6356 1290}%
\special{fp}%
\special{pa 6342 1324}%
\special{pa 6339 1332}%
\special{fp}%
\special{pa 6325 1366}%
\special{pa 6321 1373}%
\special{fp}%
\special{pa 6306 1407}%
\special{pa 6303 1414}%
\special{fp}%
\special{pa 6286 1447}%
\special{pa 6283 1455}%
\special{fp}%
\special{pa 6265 1487}%
\special{pa 6261 1494}%
\special{fp}%
\special{pa 6242 1526}%
\special{pa 6238 1533}%
\special{fp}%
\special{pa 6217 1563}%
\special{pa 6212 1570}%
\special{fp}%
\special{pa 6191 1600}%
\special{pa 6186 1607}%
\special{fp}%
\special{pa 6162 1635}%
\special{pa 6157 1641}%
\special{fp}%
\special{pa 6132 1668}%
\special{pa 6126 1674}%
\special{fp}%
\special{pa 6100 1700}%
\special{pa 6094 1706}%
\special{fp}%
\special{pa 6066 1731}%
\special{pa 6060 1736}%
\special{fp}%
\special{pa 6031 1759}%
\special{pa 6025 1764}%
\special{fp}%
\special{pa 5995 1786}%
\special{pa 5988 1790}%
\special{fp}%
\special{pa 5958 1811}%
\special{pa 5951 1815}%
\special{fp}%
\special{pa 5919 1835}%
\special{pa 5912 1839}%
\special{fp}%
\special{pa 5880 1857}%
\special{pa 5873 1860}%
\special{fp}%
\special{pa 5840 1878}%
\special{pa 5833 1882}%
\special{fp}%
\special{pa 5800 1898}%
\special{pa 5793 1902}%
\special{fp}%
\special{pa 5759 1917}%
\special{pa 5752 1921}%
\special{fp}%
\special{pa 5718 1936}%
\special{pa 5711 1939}%
\special{fp}%
\special{pa 5676 1953}%
\special{pa 5669 1956}%
\special{fp}%
\special{pa 5635 1970}%
\special{pa 5627 1973}%
\special{fp}%
\special{pa 5593 1987}%
\special{pa 5585 1989}%
\special{fp}%
\special{pa 5550 2003}%
\special{pa 5543 2005}%
\special{fp}%
\special{pa 5508 2018}%
\special{pa 5500 2021}%
\special{fp}%
\special{pa 5466 2033}%
\special{pa 5458 2036}%
\special{fp}%
\special{pa 5423 2048}%
\special{pa 5415 2051}%
\special{fp}%
\special{pa 5380 2063}%
\special{pa 5373 2065}%
\special{fp}%
\special{pa 5338 2077}%
\special{pa 5330 2080}%
\special{fp}%
% SPLINE 2 2 3 0 Black White
% 4 8090 2220 7340 2650 7010 3470 7010 3470
% 
\special{pn 8}%
\special{pn 8}%
\special{pa 8090 2220}%
\special{pa 8083 2223}%
\special{fp}%
\special{pa 8048 2237}%
\special{pa 8041 2240}%
\special{fp}%
\special{pa 8006 2253}%
\special{pa 7999 2256}%
\special{fp}%
\special{pa 7964 2270}%
\special{pa 7957 2273}%
\special{fp}%
\special{pa 7923 2287}%
\special{pa 7915 2290}%
\special{fp}%
\special{pa 7881 2304}%
\special{pa 7874 2307}%
\special{fp}%
\special{pa 7840 2322}%
\special{pa 7832 2325}%
\special{fp}%
\special{pa 7799 2341}%
\special{pa 7791 2344}%
\special{fp}%
\special{pa 7758 2360}%
\special{pa 7751 2363}%
\special{fp}%
\special{pa 7717 2379}%
\special{pa 7710 2383}%
\special{fp}%
\special{pa 7677 2400}%
\special{pa 7670 2403}%
\special{fp}%
\special{pa 7637 2421}%
\special{pa 7630 2425}%
\special{fp}%
\special{pa 7598 2443}%
\special{pa 7591 2448}%
\special{fp}%
\special{pa 7560 2467}%
\special{pa 7553 2471}%
\special{fp}%
\special{pa 7522 2492}%
\special{pa 7516 2496}%
\special{fp}%
\special{pa 7486 2518}%
\special{pa 7479 2523}%
\special{fp}%
\special{pa 7450 2545}%
\special{pa 7444 2551}%
\special{fp}%
\special{pa 7416 2575}%
\special{pa 7410 2580}%
\special{fp}%
\special{pa 7383 2605}%
\special{pa 7377 2611}%
\special{fp}%
\special{pa 7351 2638}%
\special{pa 7346 2644}%
\special{fp}%
\special{pa 7322 2672}%
\special{pa 7317 2678}%
\special{fp}%
\special{pa 7294 2708}%
\special{pa 7289 2714}%
\special{fp}%
\special{pa 7268 2744}%
\special{pa 7264 2751}%
\special{fp}%
\special{pa 7245 2783}%
\special{pa 7241 2790}%
\special{fp}%
\special{pa 7222 2822}%
\special{pa 7218 2829}%
\special{fp}%
\special{pa 7202 2862}%
\special{pa 7198 2869}%
\special{fp}%
\special{pa 7182 2903}%
\special{pa 7179 2910}%
\special{fp}%
\special{pa 7164 2944}%
\special{pa 7161 2951}%
\special{fp}%
\special{pa 7148 2986}%
\special{pa 7145 2993}%
\special{fp}%
\special{pa 7132 3028}%
\special{pa 7130 3036}%
\special{fp}%
\special{pa 7118 3071}%
\special{pa 7115 3078}%
\special{fp}%
\special{pa 7104 3114}%
\special{pa 7102 3121}%
\special{fp}%
\special{pa 7090 3157}%
\special{pa 7088 3164}%
\special{fp}%
\special{pa 7078 3200}%
\special{pa 7076 3208}%
\special{fp}%
\special{pa 7066 3244}%
\special{pa 7064 3251}%
\special{fp}%
\special{pa 7055 3287}%
\special{pa 7053 3295}%
\special{fp}%
\special{pa 7044 3331}%
\special{pa 7042 3339}%
\special{fp}%
\special{pa 7033 3375}%
\special{pa 7031 3382}%
\special{fp}%
\special{pa 7022 3418}%
\special{pa 7021 3426}%
\special{fp}%
\special{pa 7012 3462}%
\special{pa 7010 3470}%
\special{fp}%
% SPLINE 2 2 3 0 Black White
% 4 5320 2240 6130 2560 6520 3490 6520 3490
% 
\special{pn 8}%
\special{pn 8}%
\special{pa 5320 2240}%
\special{pa 5328 2242}%
\special{fp}%
\special{pa 5364 2250}%
\special{pa 5372 2252}%
\special{fp}%
\special{pa 5408 2260}%
\special{pa 5416 2262}%
\special{fp}%
\special{pa 5452 2270}%
\special{pa 5460 2272}%
\special{fp}%
\special{pa 5496 2281}%
\special{pa 5504 2283}%
\special{fp}%
\special{pa 5540 2291}%
\special{pa 5548 2293}%
\special{fp}%
\special{pa 5584 2303}%
\special{pa 5591 2305}%
\special{fp}%
\special{pa 5627 2315}%
\special{pa 5635 2317}%
\special{fp}%
\special{pa 5671 2327}%
\special{pa 5678 2329}%
\special{fp}%
\special{pa 5714 2340}%
\special{pa 5721 2343}%
\special{fp}%
\special{pa 5756 2355}%
\special{pa 5764 2358}%
\special{fp}%
\special{pa 5799 2370}%
\special{pa 5807 2373}%
\special{fp}%
\special{pa 5841 2386}%
\special{pa 5849 2389}%
\special{fp}%
\special{pa 5883 2404}%
\special{pa 5890 2407}%
\special{fp}%
\special{pa 5924 2422}%
\special{pa 5931 2426}%
\special{fp}%
\special{pa 5964 2444}%
\special{pa 5971 2447}%
\special{fp}%
\special{pa 6003 2466}%
\special{pa 6010 2470}%
\special{fp}%
\special{pa 6041 2490}%
\special{pa 6048 2494}%
\special{fp}%
\special{pa 6078 2516}%
\special{pa 6084 2521}%
\special{fp}%
\special{pa 6113 2545}%
\special{pa 6119 2550}%
\special{fp}%
\special{pa 6146 2575}%
\special{pa 6152 2581}%
\special{fp}%
\special{pa 6177 2608}%
\special{pa 6183 2614}%
\special{fp}%
\special{pa 6206 2643}%
\special{pa 6211 2649}%
\special{fp}%
\special{pa 6234 2679}%
\special{pa 6238 2685}%
\special{fp}%
\special{pa 6259 2716}%
\special{pa 6263 2723}%
\special{fp}%
\special{pa 6282 2755}%
\special{pa 6287 2761}%
\special{fp}%
\special{pa 6304 2794}%
\special{pa 6308 2801}%
\special{fp}%
\special{pa 6324 2835}%
\special{pa 6328 2842}%
\special{fp}%
\special{pa 6343 2876}%
\special{pa 6346 2883}%
\special{fp}%
\special{pa 6360 2917}%
\special{pa 6363 2925}%
\special{fp}%
\special{pa 6376 2960}%
\special{pa 6379 2967}%
\special{fp}%
\special{pa 6392 3002}%
\special{pa 6394 3010}%
\special{fp}%
\special{pa 6406 3045}%
\special{pa 6408 3053}%
\special{fp}%
\special{pa 6420 3088}%
\special{pa 6422 3096}%
\special{fp}%
\special{pa 6432 3131}%
\special{pa 6434 3139}%
\special{fp}%
\special{pa 6445 3175}%
\special{pa 6447 3182}%
\special{fp}%
\special{pa 6456 3218}%
\special{pa 6458 3226}%
\special{fp}%
\special{pa 6467 3262}%
\special{pa 6469 3270}%
\special{fp}%
\special{pa 6478 3306}%
\special{pa 6480 3314}%
\special{fp}%
\special{pa 6488 3350}%
\special{pa 6490 3358}%
\special{fp}%
\special{pa 6498 3394}%
\special{pa 6500 3402}%
\special{fp}%
\special{pa 6508 3438}%
\special{pa 6510 3446}%
\special{fp}%
\special{pa 6518 3482}%
\special{pa 6520 3490}%
\special{fp}%
% STR 2 0 3 0 Black White
% 4 6930 1390 6930 1490 2 0 0 0
% $J_1^-$
\put(69.3000,-14.9000){\makebox(0,0)[lb]{$J_1^-$}}%
% STR 2 0 3 0 Black White
% 4 6360 1380 6360 1480 2 0 0 0
% $J_2^-$
\put(63.6000,-14.8000){\makebox(0,0)[lb]{$J_2^-$}}%
% STR 2 0 3 0 Black White
% 4 7020 2680 7020 2780 2 0 0 0
% $J_2^+$
\put(70.2000,-27.8000){\makebox(0,0)[lb]{$J_2^+$}}%
% STR 2 0 3 0 Black White
% 4 6390 2680 6390 2780 2 0 0 0
% $J_1^+$
\put(63.9000,-27.8000){\makebox(0,0)[lb]{$J_1^+$}}%
% SPLINE 2 2 3 0 Black White
% 4 7990 1640 7670 2110 8040 2570 8040 2570
% 
\special{pn 8}%
\special{pn 8}%
\special{pa 7990 1640}%
\special{pa 7984 1646}%
\special{fp}%
\special{pa 7958 1672}%
\special{pa 7952 1678}%
\special{fp}%
\special{pa 7926 1704}%
\special{pa 7920 1710}%
\special{fp}%
\special{pa 7893 1736}%
\special{pa 7888 1742}%
\special{fp}%
\special{pa 7862 1769}%
\special{pa 7856 1775}%
\special{fp}%
\special{pa 7831 1802}%
\special{pa 7826 1808}%
\special{fp}%
\special{pa 7801 1836}%
\special{pa 7796 1842}%
\special{fp}%
\special{pa 7772 1871}%
\special{pa 7767 1878}%
\special{fp}%
\special{pa 7745 1908}%
\special{pa 7741 1915}%
\special{fp}%
\special{pa 7721 1946}%
\special{pa 7717 1953}%
\special{fp}%
\special{pa 7699 1986}%
\special{pa 7696 1993}%
\special{fp}%
\special{pa 7683 2028}%
\special{pa 7681 2036}%
\special{fp}%
\special{pa 7673 2073}%
\special{pa 7671 2080}%
\special{fp}%
\special{pa 7670 2118}%
\special{pa 7671 2126}%
\special{fp}%
\special{pa 7678 2162}%
\special{pa 7680 2170}%
\special{fp}%
\special{pa 7693 2205}%
\special{pa 7697 2212}%
\special{fp}%
\special{pa 7714 2245}%
\special{pa 7718 2252}%
\special{fp}%
\special{pa 7739 2283}%
\special{pa 7744 2290}%
\special{fp}%
\special{pa 7767 2319}%
\special{pa 7772 2325}%
\special{fp}%
\special{pa 7797 2353}%
\special{pa 7803 2359}%
\special{fp}%
\special{pa 7829 2385}%
\special{pa 7835 2391}%
\special{fp}%
\special{pa 7862 2417}%
\special{pa 7868 2422}%
\special{fp}%
\special{pa 7895 2447}%
\special{pa 7901 2453}%
\special{fp}%
\special{pa 7930 2477}%
\special{pa 7936 2483}%
\special{fp}%
\special{pa 7964 2507}%
\special{pa 7970 2512}%
\special{fp}%
\special{pa 7999 2536}%
\special{pa 8005 2541}%
\special{fp}%
\special{pa 8034 2565}%
\special{pa 8040 2570}%
\special{fp}%
% SPLINE 2 2 3 0 Black White
% 4 6120 900 6720 980 7360 900 7360 900
% 
\special{pn 8}%
\special{pn 8}%
\special{pa 6120 900}%
\special{pa 6128 902}%
\special{fp}%
\special{pa 6165 909}%
\special{pa 6173 910}%
\special{fp}%
\special{pa 6211 917}%
\special{pa 6219 919}%
\special{fp}%
\special{pa 6256 926}%
\special{pa 6264 928}%
\special{fp}%
\special{pa 6301 934}%
\special{pa 6309 936}%
\special{fp}%
\special{pa 6346 943}%
\special{pa 6354 944}%
\special{fp}%
\special{pa 6392 950}%
\special{pa 6400 951}%
\special{fp}%
\special{pa 6437 957}%
\special{pa 6446 958}%
\special{fp}%
\special{pa 6483 963}%
\special{pa 6491 964}%
\special{fp}%
\special{pa 6529 969}%
\special{pa 6537 969}%
\special{fp}%
\special{pa 6575 973}%
\special{pa 6583 973}%
\special{fp}%
\special{pa 6621 976}%
\special{pa 6629 977}%
\special{fp}%
\special{pa 6667 978}%
\special{pa 6675 979}%
\special{fp}%
\special{pa 6713 980}%
\special{pa 6721 980}%
\special{fp}%
\special{pa 6759 980}%
\special{pa 6767 980}%
\special{fp}%
\special{pa 6805 978}%
\special{pa 6813 978}%
\special{fp}%
\special{pa 6851 976}%
\special{pa 6859 975}%
\special{fp}%
\special{pa 6897 973}%
\special{pa 6905 972}%
\special{fp}%
\special{pa 6943 968}%
\special{pa 6951 967}%
\special{fp}%
\special{pa 6988 963}%
\special{pa 6997 961}%
\special{fp}%
\special{pa 7034 957}%
\special{pa 7042 956}%
\special{fp}%
\special{pa 7080 950}%
\special{pa 7088 949}%
\special{fp}%
\special{pa 7125 943}%
\special{pa 7133 941}%
\special{fp}%
\special{pa 7171 935}%
\special{pa 7179 934}%
\special{fp}%
\special{pa 7216 927}%
\special{pa 7224 925}%
\special{fp}%
\special{pa 7262 919}%
\special{pa 7270 918}%
\special{fp}%
\special{pa 7307 911}%
\special{pa 7315 909}%
\special{fp}%
\special{pa 7352 901}%
\special{pa 7360 900}%
\special{fp}%
% SPLINE 2 2 3 0 Black White
% 4 5370 1740 5640 2170 5410 2570 5410 2570
% 
\special{pn 8}%
\special{pn 8}%
\special{pa 5370 1740}%
\special{pa 5376 1746}%
\special{fp}%
\special{pa 5402 1774}%
\special{pa 5408 1780}%
\special{fp}%
\special{pa 5434 1808}%
\special{pa 5440 1814}%
\special{fp}%
\special{pa 5465 1843}%
\special{pa 5471 1849}%
\special{fp}%
\special{pa 5496 1878}%
\special{pa 5502 1884}%
\special{fp}%
\special{pa 5526 1914}%
\special{pa 5531 1920}%
\special{fp}%
\special{pa 5554 1951}%
\special{pa 5559 1958}%
\special{fp}%
\special{pa 5580 1990}%
\special{pa 5584 1997}%
\special{fp}%
\special{pa 5603 2030}%
\special{pa 5607 2038}%
\special{fp}%
\special{pa 5622 2073}%
\special{pa 5625 2081}%
\special{fp}%
\special{pa 5635 2118}%
\special{pa 5636 2126}%
\special{fp}%
\special{pa 5640 2165}%
\special{pa 5640 2173}%
\special{fp}%
\special{pa 5636 2211}%
\special{pa 5634 2219}%
\special{fp}%
\special{pa 5624 2256}%
\special{pa 5621 2264}%
\special{fp}%
\special{pa 5607 2300}%
\special{pa 5603 2307}%
\special{fp}%
\special{pa 5585 2341}%
\special{pa 5581 2348}%
\special{fp}%
\special{pa 5560 2380}%
\special{pa 5555 2387}%
\special{fp}%
\special{pa 5533 2419}%
\special{pa 5528 2425}%
\special{fp}%
\special{pa 5505 2456}%
\special{pa 5500 2462}%
\special{fp}%
\special{pa 5475 2492}%
\special{pa 5470 2498}%
\special{fp}%
\special{pa 5446 2528}%
\special{pa 5440 2535}%
\special{fp}%
\special{pa 5415 2564}%
\special{pa 5410 2570}%
\special{fp}%
% STR 2 0 3 0 Black White
% 4 6650 2830 6650 2930 2 0 0 0
% $I_1$
\put(66.5000,-29.3000){\makebox(0,0)[lb]{$I_1$}}%
% STR 2 0 3 0 Black White
% 4 7490 2060 7490 2160 2 0 0 0
% $I_2$
\put(74.9000,-21.6000){\makebox(0,0)[lb]{$I_2$}}%
% STR 2 0 3 0 Black White
% 4 6610 1010 6610 1110 2 0 0 0
% $I_3$
\put(66.1000,-11.1000){\makebox(0,0)[lb]{$I_3$}}%
% STR 2 0 3 0 Black White
% 4 5700 2120 5700 2220 2 0 0 0
% $I_4$
\put(57.0000,-22.2000){\makebox(0,0)[lb]{$I_4$}}%
\end{picture}}%
\caption{$\rho(\gamma_1)=\rho(\beta^2\alpha\beta\alpha)$ maps $\alpha\beta^2\alpha$
to $\beta^2$ and $\alpha\beta^2$ to
 $\beta^2\alpha$. $\rho(\gamma_2)=\rho(\alpha\beta\alpha\beta^2)$ maps $\beta$ to
$\alpha\beta\alpha$ and $\beta\alpha$ to $\alpha\beta$.
}
\end{figure}
Define  open
intervals $J_1^-$, $J_1^+$, $J_2^-$ and $J_2^+$, slightly bigger than
$K_1^-$, $K_1^+$, $K_2^-$ and $K_2^+$. See Figure 4. ($J_i^\pm$ and
$K_i^\pm$ are actually intervals of $S^1$.)
Recall the interval $I_1$, the gap of $\MM$ containing $y_0$, and let
$I_2=\rho(\alpha\beta^2)I_1$, $I_3=\rho(\alpha)I_1$ and 
$I_4=\rho(\beta\alpha) I_1$ (Figure 4). They are gaps of the minimal set
$\MM$ and the orbit $\rho(\oG)y_0$ is discrete
in each of them. 
{ 
Since the stabilizer of $I_1$ is generated by $\alpha\beta$, the three
points in $I_1$, $\beta^2\alpha$, $e$ and $\alpha\beta$ are consecutive 
points of $\rho(\oG)y_0$ contained in $I_1$.
 Their images by $\rho(\alpha)$, $\alpha\beta^2\alpha$, $\alpha$ and $\beta$
are consecutive points of $\rho(\oG)y_0$ contained in $I_3$.
Likewise $\alpha\beta\alpha$ and $\alpha\beta^2$ are consecutive in
$I_2$,
and $\beta\alpha$ and $\beta^2$ are consecutive in $I_4$.
One chooses $J_2^+$ so that the point $\alpha\beta$ (resp.
 $\alpha\beta\alpha$) is the leftmost (resp.\ rightmost) point of 
$\rho(\oG)y_0\cap J_2^+$.} More generally, the intervals
$J_i^\pm$ are so chosen that points of $\partial K_i^{\pm}$ are
extremal in $\rho(\oG)y_0\cap J_i^{\pm}$.
These points are called the {\em guardians} of the interval with respect to $\rho$.
Notice that there are just two
points of $\rho(\oG)y_0$ outside of $\bigcup_{i,\pm}J_i^\pm$, namely
$y_0$ and $\rho(\alpha)y_0$ (denoted by $e$ and $\alpha$ in Figure 4).

Define a neighbourhood $U$ of $\rho$
in ${\rm Hom}(\oG,{\rm  Homeo}_+(S^1))$ such that each element $\rho'\in U$
satisfies the following conditions.
\\
(1) $\rho'(\gamma_i)$ maps the closed set $S^1\setminus J_i^-$
 into the open set $J_i^+$ ($i=1,2$).
\\
(2) The configuration of $\rho'(g)y_0$ in  $S^1$ for ten elements  
\begin{equation}\label{lasteq}
g=e,\alpha\beta,\alpha\beta\alpha,\alpha\beta^2,\alpha\beta^2\alpha,\alpha,\beta,\beta\alpha,\beta^2,\beta^2\alpha
\end{equation}
is the same as for $\rho$ (Figure 4), as well as their configuration with respect to the four
intervals $J_i^\pm$.

Take a neighbourhood $V$ of $c$ as in Proposition \ref{continuous}, 
and for any $c'\in V$, let $\rho'$ be a conjugate of a dynamical
realization of $c'$ which is contained in $U$.
Then by the ping-pong argument,
 the circular order of the orbit $\rho'(\oG)y_0$ is uniquely determined.
Let us show this a bit in detail. Call ten points $\rho'(g)y_0$ ($g$ as in (\ref{lasteq}))
{\em of the first generation}.
The images of points of the first generation by $\gamma_i^{\pm1}$ which
are not themselves of first generation are called {\em of second generation}.
Then the configuration of the points of first and second generations are
uniquely determined. In fact, the guardians $\rho'(\alpha\beta^2)y_0$
and $\rho'(\alpha\beta^2\alpha)y_0$ of the interval $J_1^-$ is mapped by 
$\rho'(\gamma_1)$ to the guardians $\rho'(\beta^2\alpha)y_0$ and
 $\rho'(\beta^2)y_0$ of $J_1^+$, and all the other eight points are mapped into
the interval in $J_1^+$ bounded by the latter guardians. The same is true for
$\gamma_1^{-1}$ and $\gamma_2^{\pm1}$. Since $\gamma_i^{\pm1}$ are
orientation preserving, the configuration of the points of the first and second
generation is uniquely determined.
Next we define points of third generation in a similar way.
These points are contained in $\bigcup_{i,\pm}J_i^\pm$. For example,
those contained in $J_1^+$ are the images
of the points of second generation in $S^1\setminus J_1^-$ by
$\rho'(\gamma_1)$. The configuration of these points, together with
the points
of first and second generations, is uniquely determined. 
{ 
On the other hand, any point in $\rho'(\oG)y_0$ is a point of some
generation. This is because the ten elements of (9.1)
exhaust the set of right cosets $[\oG,\oG]\setminus \oG$, and the elements
$\gamma_1$ and $\gamma_2$ generate $[\oG,\oG]$.
}
Therefore
continuing this way, we see that the
natural circular order of the whole orbit $\rho'(\oG)y_0$ is uniquely determined,
that is, the same as $\rho(\oG)y_0$.
 This shows that $c'=c$, i.e, $c$ is isolated.
Define $c^{(1)}$ in Theorem \ref{t92} to be this $c$.

\bigskip
For $k>1$, denote by $p_k:S^1\to S^1$ the $k$-fold covering
map. A  representation $\rho^{(k)}:\overline G\to{\rm Homeo}_+(S^1)$ is called
a {\em $k$-fold lift} of $\rho$ if $p_k\rho^{(k)}(g)=\rho(g)p_k$ holds for
 any $g\in \overline G$.
There is a $k$-fold lift $\rho^{(k)}$ of our representation $\rho$ if and only if $k\equiv \pm1$ mod 6,
and it is unique if it exists. 
{ Computation shows that if $k=6\ell\pm 1$, then 
\begin{equation}\label{rot}
 {\rm rot}(\rho^{(k)}(\alpha\beta))=\mp \ell/k. 
\end{equation}}
Notice that $(k,\ell)=1$.
We fix such $k$.

Let $y_0^{\mu}$, $1\leq \mu\leq k$,
be the lifts of the point $y_0$ by $p_k$.
The the natural circular order of the orbit $\rho^{(k)}(\oG)y_0^\mu$ in
$S^1$ is
the same for any $\mu$.
Denote  it by $c^{(k)}\in CO(\oG)$. 
{  It is easy to show that $\rho^{(k)}$ is tight at $y_0^\mu$.
Thus $\rho^{(k)}$ is topologically conjugate to a dynamical realization
of $c^{(k)}$ at $y_0^\mu$ by an orientation and base point preserving homeomorphism.}
Let us show that $c^{(k)}$ is isolated.
Let 
$J_{i,\mu}^{\pm}$ ($\mu=1,\ldots,k$) be the connected components of 
$p_k^{-1}(J_i^{\pm})$.

Define a neighbourhood $U^{(k)}$ of $\rho^{(k)}$
such that each element $\rho'\in U^{(k)}$
satisfies the following conditions.
\\
(1) $\rho'(\gamma_i)$ maps each component of $S^1\setminus \bigcup_\mu J_{i,\mu}^-$
 into $J_{i,\nu}^+$ ($i=1,2$), where $\nu$ is determined so that
$\rho^{(k)}(\gamma_i)$ maps the same component
 into $J_{i,\nu}^+$.
\\
(2) The configuration of $10k$ points
$\rho'(g)y_0^\mu$ in  $S^1$ ($g$  as in (\ref{lasteq}) and  $1\leq\mu\leq k$)
is the same as $\rho^{(k)}$.
Their configuration relative to $J^\pm_{i,\mu}$ is also the same.

Then the same ping-pong argument as before shows that the natural circular order of
$\rho'(\oG)(p_k^{-1}(y_0))$ for $\rho'\in U^{(k)}$ is uniquely determined.
In particular, the natural circular order of $\rho'(\oG)y_0^\mu$ is 
the same as for $\rho^{(k)}$, showing that $c^{(k)}$ is
isolated.

\bigskip
Finally let us show that $c^{(k)}$'s are not the automorphic images of
the others, by considering their linear parts. In $\oG=PSL(2,\Z)$,
any element of infinite order is a multiple of a unique primitive
element. This can be shown by considering the modular representation
$\rho_M$: the fixed point set of any element of infinite order is
 either one point of 
$\partial\mathbb H$ or a two point set of $\partial\mathbb H$, and the
isotropy group of the fixed point set is infinite cyclic.

As we have seen above, the linear part of $c^{(1)}$ is generated by a
primitive element 
$(\alpha\beta)^{\pm1}$. The equality (\ref{rot}) 
shows that the linear part of $c^{(k)}$ is generated
by $(\alpha\beta)^{\pm k}$. For different choices of $k$ and $k'$, there
is no automorphism of $\oG$ which maps $(\alpha\beta)^{\pm k}$ to
$(\alpha\beta)^{\pm k'}$. 
{  In fact, any automorphism of $\oG$ maps
a primitive element to a primitive element, and
thus it should map $(\alpha\beta)^{\pm k}$ to a $k$-multiple of a
primitive element.}
This finishes the proof of Theorem \ref{t92}.

\begin{remark}
 The left order $\lambda=\pi^*c^{(1)}\in LO(G)$ is
the Dubrovina-Dubrovin order, since it satisfies
$e<_\lambda a<_\lambda b$. It can be shown that $\lambda'=\pi^*c^{(5)}$
is the unique left order which satisfies
$$
(ab)^5t^{-4}<_{\lambda'}e<_{\lambda'}a<_{\lambda'}(ab)^5t^{-4}a.$$
\end{remark}

\end{document}